\documentclass[preprint,12pt]{elsarticle}

\usepackage{amssymb}
\usepackage{hyperref}
\usepackage{setspace}
\usepackage{amsmath}
\usepackage{microtype}
\usepackage{tikz}
\usetikzlibrary{arrows}
\usetikzlibrary{arrows.meta}
\usepackage{amssymb}
\usepackage{multirow}
\usepackage{pgfplots}
\usepackage{pdfpages}
\usepackage{xcolor}
\usepackage{pdflscape}
\pgfplotsset{compat=1.4}
\usetikzlibrary{patterns}

\newtheorem{Remark}{\bf Remark}

\journal{arXiv}

\begin{document}
\begin{frontmatter}



\title{Shock-fronted travelling waves in a reaction-diffusion model with nonlinear forward-backward diffusion}


\author{Yifei Li$^1$, Peter van Heijster$^{1,2}$,  Matthew J. Simpson$^1$ and Martin Wechselberger$^3$}

\address{$^{1}$School of Mathematical Sciences, Queensland University of Technology, Brisbane, QLD 4001,
Australia\\
$^{2}$Biometris, Wageningen University and Research, Wageningen, Netherlands\\
$^{3}$School of Mathematics and Statistics
University of Sydney, Sydney, NSW 2006,
Australia}

\begin{abstract}
Reaction-diffusion equations (RDEs) are often derived as continuum limits of lattice-based discrete models. Recently, a discrete model which allows the rates of movement, proliferation and death to depend upon whether the agents are isolated has been proposed, and this approach gives various RDEs where the diffusion term is convex and can become negative (Johnston et al., Sci. Rep. 7, 2017), i.e. forward-backward diffusion. Numerical simulations suggest these RDEs support shock-fronted travelling waves when the reaction term includes an Allee effect. In this work we formalise these preliminary numerical observations by analysing the shock-fronted travelling waves through embedding the RDE into a larger class of higher order partial differential equations (PDEs). Subsequently, we use geometric singular perturbation theory to study this larger class of equations and prove the existence of these shock-fronted travelling waves. Most notable, we show that different embeddings yield shock-fronted travelling waves with different properties.
\end{abstract}




\end{frontmatter}

\section{Introduction}
\begin{figure}
\centering
\hspace*{-1cm}
\includegraphics[width=0.8\textwidth]{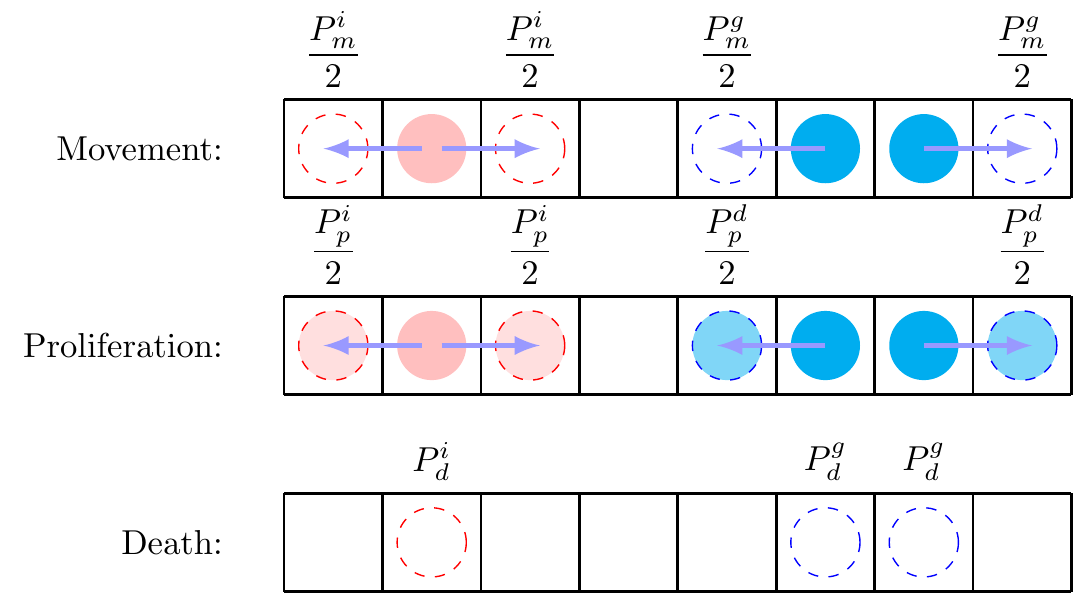}
\caption{Schematic depiction of the evolution rules of the lattice-based model introduced in \cite{johnston2017co}. Pink discs represent isolated agents and blue discs represent grouped agents. During each time step of duration $\tau$, isolated agents attempt to move to nearest neighbour lattice sites with a probability $P^i_m$, to proliferate to form new agents in neighbour sites with a probability $P^i_p$ and to die with a probability $P^i_d$. Similarly, grouped agents attempt to move to neighbour sites with a probability $P^g_m$, to proliferate to form new agents in neighbour sites with a probability $P^g_p$ and to die with a probability $P^g_d$. The attempts that would move to an occupied site or place an agent on an occupied site are aborted.}
\label{LatticePicture}
\end{figure}

Reaction-diffusion equations
(RDEs) are widely used to study population dynamics in cell biology and ecology
\cite{murray2002mathematical}. Often, $U(x,t)$ represents a population density and provides a macroscopic description of individual behaviour. For RDEs established from the continuum limit of stochastic models, a solution of the RDE not only shows the macroscopic evolution of $U(x,t)$, but it also reflects how microscopic behaviour of individuals influences the macroscopic outcomes~\cite{anguige2009one,deroulers2009,johnston2017co,johnston2012mean,2019arXiv190310090L,Simpson2010invasion}. 
Johnston et~al.~\cite{johnston2017co}
 introduced a lattice-based stochastic model
to study how a population of individuals can undergo motility, proliferation and death events with the aim of studying biological and ecological invasion, see Figure \ref{LatticePicture}.
By considering different behaviours of isolated and grouped agents, including motility, proliferation and death events, an RDE with a nonlinear diffusivity function and a logistic or Allee type reaction term was derived as the continuum limit. 
In particular,
\begin{equation}
\label{RDE_biology1}
    \begin{aligned}
    \frac{\partial U}{\partial t}=\frac{\partial}{\partial x}\left(D(U)\frac{\partial U}{\partial x}\right)+{R}\left(U\right),
    \end{aligned}
\end{equation} 
where $U(x,t)$ represents the total population density at position $x\in\mathbb{R}$ and time $t\in\mathbb{R_+}$. 

The nonlinear diffusivity function is given by
\begin{equation}
\label{D(u)}
D\left(U\right)=3(D_i-D_g)U^2-4(D_i-D_g)U+D_i,
\end{equation}
where $D_i\ge0$ and $D_g\ge0$ are diffusivities of the isolated and grouped agents, respectively. When $D_i>4D_g$, $D(U)$ has two real roots, $\alpha$ and $\beta$, which are centred around $2/3$, and are given by
\begin{align}
    \label{AB}
        \alpha=\frac{2}{3}\left(1-\sqrt{\frac{D_i-4D_g}{4(D_i-D_g)}}\right),\quad
        \beta=\frac{2}{3}\left(1+\sqrt{\frac{D_i-4D_g}{4(D_i-D_g)}}\right),
\end{align}
and $D(U)<0$ for $U \in (\alpha,\beta)$. 
While the negativity of a nonlinear diffusivity function is sometimes related to aggregation in the underlying discrete model~\cite{Simpson2010motile}, here it is actually a macroscopic effect of the isolated and the grouped motility of the agents, together with competition for space, that leads to a net \emph{aggregation effect}~\cite{johnston2017co}.
The condition $D_i>4D_g$ implies that the motility rate of isolated agents is greater than the motility rate of grouped agents, which is consistent with the common biological observation that isolated \textit{leader cells} are more motile than \textit{follower cells} \cite{Poujade15988,Simpson2014Plos}. Note that $D_i$ and $D_g$ are related to $P_m^i$ and $P_m^g$, respectively, in the lattice-based model in Figure \ref{LatticePicture}. Full details of the discrete model and the continuum limit derivation are given in~\cite{johnston2017co}. 

\begin{figure}
\centering
\includegraphics[width=0.9\textwidth]{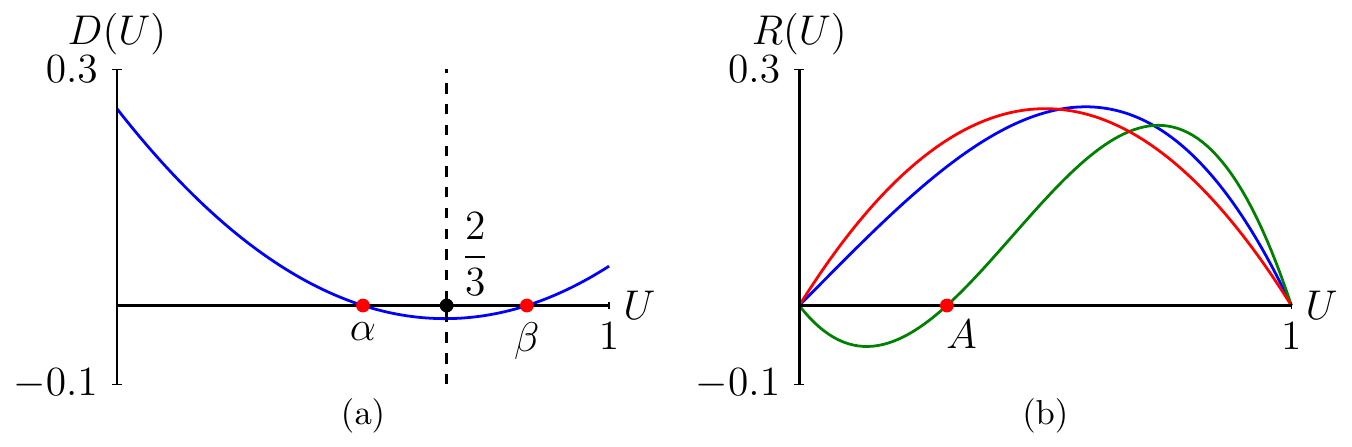}
\caption{(a) The nonlinear diffusivity function $D(U)$ (\ref{D(u)}) centred around $2/3$ (dashed line). (b) The reaction term $R(u)$ corresponding to the logistic growth with $\lambda_g=1$ (red), the weak Allee effect with $r=0.8$ and $A=-0.9$ (blue), and the strong Allee effect with $r=3$ and $A=0.3$ (green).} 
\label{TypesofAlleeeffect}
\end{figure}

The reaction term, whose parameters are also related to parameters in the lattice-based model depicted in Figure~\ref{LatticePicture}, is given by 
\begin{equation}
\label{originalR(u)}
R\left(U\right)=\lambda_gU(1-U)+(\lambda_i-\lambda_g-K_i+K_g)U(1-U)^2-K_gU,
\end{equation}
where $\lambda_i\ge0$ and $\lambda_g\ge0$ are the proliferation rates of isolated and grouped agents, respectively; $K_i\ge0$ and $K_g\ge0$ are the death rates of isolated and grouped agents, respectively \cite{johnston2017co}. If the proliferation mechanism is the same for isolated and grouped agents and no death event occurs, that is, $\lambda_i=\lambda_g$ and $K_i=K_g=0$, then (\ref{originalR(u)}) simplifies to a logistic reaction term
\begin{equation}
\label{logisticR(u)}
    R(U)=\lambda_gU(1-U).
\end{equation}
If the proliferation and death mechanisms are either competitive or co-operative, that is, $\lambda_i\ne\lambda_g$ and $K_i\ne K_g$ \cite{stephens1999allee,taylor2005}, then the reaction term takes the form of an Allee effect \cite{taylor2005}. For simplicity, but without loss of generality, we assume $K_g=0$.\footnote{Although we assume $K_i=K_g=0$ to obtain the logistic reaction term (\ref{logisticR(u)}) and $K_g=0$ to obtain the Allee reaction term (\ref{AlleeR(u)}), similar reaction terms are obtained without these assumptions by scaling the population density $U(x,t)$ \cite{johnston2017co}.} Subsequently, (\ref{originalR(u)})
simplifies to
\begin{equation}
\label{AlleeR(u)}
R(U)=rU(1-U)(U-A),
\end{equation}
where $r=K_i-\lambda_i+\lambda_g$ is the intrinsic growth rate and $A=1-\lambda_g/r$ is the Allee parameter \cite{johnston2017co}. If $r>\lambda_g$, which is equivalent to $K_i>\lambda_i$ and thus implies that isolated agents have a higher death rate than proliferation rate, then $0<A<1$ and $R(U)<0$ in $(0,A)$ and $R(U)>0$ in $(A,1)$. This represents the strong Allee effect. Conversely, if $0<r<\lambda_g$, which implies that isolated agents have a higher birth rate than death rate, then $A<0$ and $R(U)>0$ in $(0,1)$. This is called the weak Allee effect \cite{johnston2017co}. See Figure \ref{TypesofAlleeeffect} for the different potential forms of $R(U)$. For simplicity, we assume that $A\ne \alpha$ and $A\ne \beta$. 

\begin{landscape}
\begin{figure}
\centering
\vspace*{-2cm}
\includegraphics{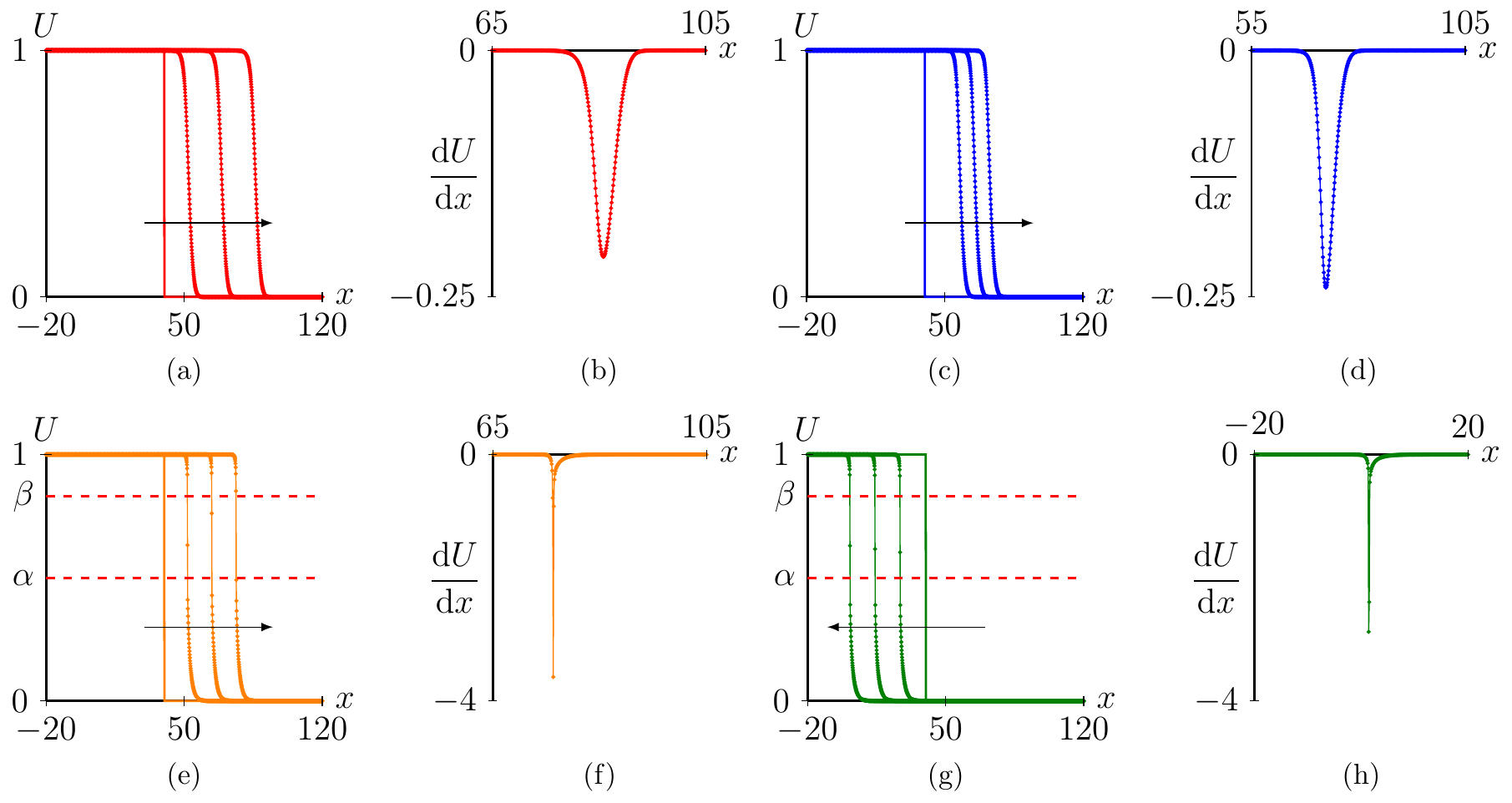}
\caption{Travelling wave solutions of (\ref{RDE_biology1}) evolving from a Heavidside initial condition, $U=1$ for $x\in[-20,40]$ and $U=0$ for $x\in(40,120]$, with $D_i=0.25$ and $D_g=0.05$ such that $D(U)=0.6U^2-0.8U+0.25$ and for different $R(U)$. 
The implicit finite difference method we used had no-flux boundary conditions, a time step $\delta t=0.01$, a space step $\delta x=0.1$ and an error tolerance $10^{-6}$. With a logistic type $R(U)=0.75U(1-U)$, (a) shows a travelling wave solution with speed $c \approx 0.86$ at $t_1=20$, $t_2=40$, $t_3=60$. With a weak Allee type $R(U)=0.5U(1-U)(U+0.2)$, where $\lambda_i=0.5$, $\lambda_g=0.6$ and $K_i=0.4$, (c) shows a smooth travelling wave solution with $c \approx 0.99$ at $t_1=15$, $t_2=30$, $t_3=45$. With a strong Allee type $R(U)=0.5U(1-U)(U-0.2)$, where $\lambda_i=0.4$, $\lambda_g=0.4$ and $K_i=0.5$, (e) shows a shock-fronted travelling wave solution with positive speed $c=0.0123$ at $t_1=1000$, $t_2=2000$, $t_3=3000$. With a different strong Allee type $R(U)=0.3U(1-U)(U-1/3)$, where $\lambda_i=0.4$, $\lambda_g=0.2$ and $K_i=0.5$, (h) shows a shock-fronted travelling wave solution with negative speed $c=-0.0127$ at $t_1=500$, $t_2=1000$, $t_3=1500$. (b) \& (d) show the derivatives of the last simulated travelling wave solutions in (a) \& (c) and highlight the smoothness of the waves. (f) \& (h) show the derivatives of the last simulated travelling wave solutions in (e) \& (g) and highlight the shocks in these waves. 
}
\label{Panel_TW1}
\end{figure}
\end{landscape}

Understanding travelling wave solutions, that is, solutions that propagate through space with a fixed shape and a constant speed, is important in the study of biological and ecological invasion processes \cite{DUCROT2021132730,ELHACHEM2020132639,harley2014existence,sewalt2016influences}. In this work, we are interested in travelling wave solutions supported by \eqref{RDE_biology1} with $D_i>4D_g$, such that the nonlinear diffusivity function is negative for $U\in(\alpha,\beta)$, see \eqref{AB}. In this case, the nonlinear diffusivity function can be written as 
\begin{equation}
\label{D(u)2}
D\left(U\right)=k (U-\alpha)(U-\beta)\,, \qquad k:=3(D_i-D_g)\,.
\end{equation}
With the implicit finite difference method introduced in \cite{johnston2017co}, numerical solutions of (\ref{RDE_biology1}) with $D(U)$ as in (\ref{D(u)2}) and with either logistic or weak Allee forms for $R(U)$ lead to smooth travelling wave solutions with positive speeds, while simulations of (\ref{RDE_biology1}) with $D(U)$ as in (\ref{D(u)2}) and strong Allee forms for $R(U)$ lead to shock-fronted travelling wave solutions with either positive or negative speeds~\cite{johnston2017co, 2019arXiv190310090L}, see Figure \ref{Panel_TW1} for different travelling wave solutions at $t=t_1, t_2, t_3$, with $t_1<t_2<t_3$. To calculate the wave speed, we locate the front of the wave by looking for the left-most coordinate $x_l$ satisfying $U(x_l,t)<10^{-3}$. 
Then, we estimate the speed from the distance the front of the wave has travelled from $t_2$ to $t_3$.
Interestingly, the speeds of the shock-fronted travelling wave solutions are much smaller than the speeds of the smooth travelling wave solutions, which potentially indicates that the mechanisms giving rise to shock-fronted 
travelling waves are fundamentally different to the mechanisms that give rise to smooth travelling waves. Note that with the nonlinear diffusivity function $D(U)$ centred around $2/3$ given by \eqref{D(u)} we only observe shock-fronted travelling wave solutions with the strong Allee effect. However, a forward-backward-forward nonlinear diffusivity function which is not centred around $2/3$ may also lead to shock-fronted travelling wave solutions with logistic or weak Allee forms of $R(U)$, see Figures 9 and 10 in \cite{2019arXiv190310090L} for an example.

Ferracuti et al. \cite{ferracuti2009travelling} showed that there exist smooth travelling wave solutions of (\ref{RDE_biology1}) with logistic $R(U)$ for a range of positive wave speeds based on the \textit{comparison method} \cite{ARONSON197833}. Kuzmin and Ruggerini \cite{Kuzmin2011819} provided necessary conditions for the existence of smooth travelling wave solutions of (\ref{RDE_biology1}) with $R(U)$ that takes the form of an Allee effect and the speed of the wave can be either negative or positive according to the shape of $D(U)$ and $R(U)$. However, to the best of our knowledge, the existence of shock-fronted travelling wave solutions to (\ref{RDE_biology1}) with $D(U)$ (\ref{D(u)2}) and $R(U)$ taking the form of an Allee effect is an open question. The methods used in \cite{Kuzmin2011819} can also be used to identify necessary conditions for the existence of shock-fronted travelling wave solutions. In particular, let $U_1$ and $U_2$ (with $U_1 \leq \alpha < \beta \leq U_2$) denote the $U$-values at 
\textls[40]{the endpoints of the shock, then a necessary condition for the existence of a} 

\noindent
monotonically decreasing shock-fronted left-travelling wave solution is
\begin{equation}
\label{necessarycondition1}
    \int_0^{U_a} D(U)R(U)\textrm{d}U<0\,, 
\end{equation}
for all $U_a\in(0,U_1)$.  Similarly, a necessary condition for the existence of a monotonically decreasing shock-fronted right-travelling wave solution is
\begin{equation}
\label{necessarycondition2}
    \int_{U_b}^1 D(U)R(U)\textrm{d}U>0\,,
\end{equation}
for all $U_b\in(U_2,1)$. We refer to \ref{A:SH} for a derivation, inspired by \cite{Kuzmin2011819}, of these necessary conditions. Heuristically this means that for shock-fronted travelling wave solutions with a positive speed, the $A$ value in $R(U)$, in order to satisfy (\ref{necessarycondition2}), should not be too close to $1$. Since $A=1-\lambda_g/r$ where $r=K_i-\lambda_i+\lambda_g$, this implies that a relatively much higher death rate of isolated agents compared to the birth rate of isolated agents will not result in a successful invasion event. Similarly, a very small birth rate of grouped agents will also not result in a successful invasion event.

In \cite{2019arXiv190310090L}, we derived, among other things, the same condition as in \cite{ferracuti2009travelling} for the existence of smooth travelling wave solutions of $(\ref{RDE_biology1})$ with logistic $R(U)$
by using a geometric approach. 
Furthermore, geometric approaches have been used to study shock-fronted travelling wave solutions. For example, in~\cite{harley2014existence,sewalt2016influences}, the authors studied shock-fronted travelling wave solutions in an advection-reaction-diffusion equation for malignant tumour invasion using geometric singular perturbation theory (GSPT) \cite{fenichel1979geometric,hek2010geometric,jones1995geometric} and canard theory~\cite{Szmolyan2001,wechselberger2012propos,Wechselberger2010}.
In this work, we use GSPT to further explore the existence of shock-fronted travelling wave solutions of (\ref{RDE_biology1}) with $\lim_{x\to-\infty}U(x,t)=1$ and $\lim_{x\to\infty}U(x,t)=0.$\footnote{$U \equiv 0$ and $U \equiv 1$ are both constant solutions of (\ref{RDE_biology1}) with the logistic $R(U)$ and the weak and strong Allee type $R(U)$. However, $U \equiv 0$ and $U \equiv 1$ are both PDE stable for the strong Allee type $R(U)$, while only $U \equiv 1$ is PDE stable for the logistic $R(U)$ and the weak Allee type $R(U)$. Therefore, it is no surprise that left-moving traveling wave solutions are only found for the strong Allee type $R(U)$, see Figure~\ref{Panel_TW1}.} Moreover, we assume $D_i>4D_g$ -- such that $D(U)<0$ for $U \in (\alpha,\beta)$ -- and 
$K_g=0$ and $r>\lambda_g$ -- such that we have a strong Allee effect type $R(U)$.
To apply GSPT, we smooth out the shock and regularise~(\ref{RDE_biology1}) by adding a small higher order perturbation term. This embeds (\ref{RDE_biology1}) into a larger class of PDEs. 
Regularisation of RDEs is typically considered in one of two ways~\cite{padra3n2004effect,pego1989front}. The first method of regularisation accounts for non-local effects by adding a small fourth-order spatial derivative term \cite{pego1989front,witelski1995shocks}. In particular, equation (\ref{RDE_biology1}) becomes
\begin{equation}
\label{RDE_biology3}
    \frac{\partial U}{\partial t}=\frac{\partial}{\partial x}\left(D(U)\frac{\partial U}{\partial x}\right)+R(U)-\varepsilon^2\frac{\partial^4U}{\partial x^4}, \quad 0\leq\varepsilon\ll1.
\end{equation}
The second method of regularisation accounts for viscous relaxation by adding a 
small mixed derivative term \cite{novick1991stable,padra3n2004effect,WitelskiInternalLayers}. In particular, equation (\ref{RDE_biology1}) becomes
\begin{equation}
\label{Another1}
    \begin{aligned}
    \frac{\partial U}{\partial t}=\frac{\partial}{\partial x}\left(D(U)\frac{\partial U}{\partial x}\right)+R(U)+\varepsilon\frac{\partial^3U}{\partial x^2\partial t}, \quad 0\leq\varepsilon\ll1.
    \end{aligned}
\end{equation}
It is important that the sign of the perturbation terms in (\ref{RDE_biology3}) and (\ref{Another1}) is such that setting $\varepsilon>0$ generally leads to well-posed problems. However, see \cite{2019arXiv190310090L} and references therein, for a further discussion related to the well-posedness of (\ref{RDE_biology1}). Also note that other types of regularisations have been used to smooth out shocks \cite{Barenblatt1993}.

In \S\ref{sec2}, we study travelling wave solutions of (\ref{RDE_biology3}) and first derive 
a higher-dimensional slow-fast system of ordinary differential equations (ODEs). The related reduced singular limit ODE systems give useful information of underlying shock-fronted travelling wave solutions of (\ref{RDE_biology1}) and (\ref{RDE_biology3}) based on GSPT and Fenichel theory \cite{fenichel1979geometric}. Because the reduced systems are algebraically intractable, we use a numerical ODE solver to determine the speed of the shock-fronted travelling wave solutions. 
In \S\ref{sec3}, we use a similar approach to establish a different higher-dimensional system of ODEs based on the viscous relaxation PDE (\ref{Another1}) and find shock-fronted travelling wave solutions with different properties. Note that in this case, GSPT has to be extended since the critical manifold loses normal hyperbolicity near a fold point.
Although (\ref{RDE_biology3}) and (\ref{Another1}) are the same in the singular limit $\varepsilon=0$, they yield shock-fronted travelling wave solutions with different speeds and different shock sizes when $\varepsilon>0$. Finally, we discuss various extensions of the current work including the relationship between the discrete model and the continuous description, the option of including different regularisation terms, the possibility of shock-fronted travelling wave solutions with logistic $R(U)$ and the spectral stability of travelling wave solutions of~(\ref{RDE_biology1}).

\begin{Remark}
In the remainder of this article we will use nonlinear diffusivity functions $D(U)$ \eqref{D(u)} and reaction terms $R(U)$ \eqref{originalR(u)} that are larger than the $D(U)$ and $R(U)$ used in Figure~\ref{Panel_TW1} to generate larger speeds. As the model based on \eqref{RDE_biology1} is dimensionless, those larger parameters in $D(U)$ and $R(U)$ still correspond to the parameters introduced in the latticed-based model in Figure~\ref{LatticePicture} upon rescaling space and/or time. However, note that the connection between the discrete and continuum models is only accurate when the rate of motility of both the grouped and isolated agents is much greater than rate of proliferation and death of both the grouped and isolated agents  and this should be kept in mind when rescaling space and/or time.
For more details, see \cite{johnston2017co}.
\end{Remark}

\section{Non-local regularisation}
\label{sec2}
In this section, we look for shock-fronted travelling wave solutions of $(\ref{RDE_biology3})$ connecting $U=1$ to $U=0$. We first introduce a travelling wave coordinate to transform (\ref{RDE_biology3}) into a fourth-order ODE. Next, we use a dynamical system approach to transform the ODE into a four-dimensional singular perturbed slow-fast system. The four-dimensional system has two equivalent forms as $\varepsilon\ne0$. However, these forms produce different lower-dimensional subsystems, called the \textit{reduced problem} and the \textit{layer problem} in the singular limit $\varepsilon=0$.\footnote{The reduced problem is called the \textit{slow reduced system} and the layer problem is called the \textit{fast reduced system} in studies of phase separation, see for example \cite{doelman2009pulse}.} 
The concatenation of solutions of each of the subsystems yields a solution of the four-dimensional system in the singular limit. We give an outline, and conclude based on GSPT, that it persists for $\varepsilon$ sufficiently small in the full four-dimensional system. This solution corresponds to a travelling wave solution of (\ref{RDE_biology3}). 
\subsection{Preliminary observations}
A travelling wave solution of (\ref{RDE_biology3}) is a solution  of the form $U(x,t)=u(x-ct)=u(z)$, where $c\in\mathbb{R}$ is the constant speed of the travelling wave solution and $z=x-ct$ is the travelling wave coordinate. Writing (\ref{RDE_biology3}) in its travelling wave coordinate leads to
\begin{equation}
\label{GSPT_ALLEE_ODE1_partial}
    \frac{\partial U}{\partial t}=\frac{\partial }{\partial z}\left(\varepsilon^2\frac{\partial^3U}{\partial z^3}-cU-\frac{\partial}{\partial z}(F(U))\right)-R(U),
\end{equation}
where $F(U)=\int D(U) \textrm{d} U$ and the reaction term $R(U)$ \eqref{AlleeR(u)} is of strong Allee effect type. A travelling wave solution $u(z)$ is a stationary solution to (\ref{GSPT_ALLEE_ODE1_partial}) that asymptotes to one as $z\to-\infty$ and to zero as $z\to\infty$. Thus, it satisfies
\begin{equation}
\label{GSPT_ALLEE_ODE1}
    0=\frac{\textrm{d}}{\textrm{d}z}\left(\varepsilon^2\frac{\textrm{d}^3u}{\textrm{d}z^3}-cu-\frac{\textrm{d}}{\textrm{d}z}(F(u))\right)-R(u).
\end{equation}
Upon defining
\begin{equation}
\label{definationtoODE1}
    p:=\varepsilon^2\frac{\textrm{d}^3u}{\textrm{d}z^3}-cu-\frac{\textrm{d}}{\textrm{d}z}(F(u)),\quad
    v:=\varepsilon^2\frac{\textrm{d}^2u}{\textrm{d}z^2}-F(u),\quad
    w:=\varepsilon\frac{\textrm{d}u}{\textrm{d}z},
\end{equation}
(\ref{GSPT_ALLEE_ODE1}) transforms into a four-dimensional singular perturbed slow-fast dynamical system
\begin{equation}
\label{GSPT_ALLEE_SYSTEMSLOW}
    \left\{\begin{aligned}
  \varepsilon &\frac{\textrm{d}u}{\textrm{d}z}&&=w,\\
  \varepsilon &\frac{\textrm{d}w}{\textrm{d}z}&&=v+F(u),\\
  &\frac{\textrm{d}p}{\textrm{d}z}&&=R(u),\\
  &\frac{\textrm{d}v}{\textrm{d}z}&&=p+cu.
\end{aligned}\right.
\end{equation} 
Here, $(u,w)\in\mathbb{R}^2$ are fast variables and $(p,v)\in\mathbb{R}^2$ are slow variables. By using a stretched, or fast variable, $\xi=z/\varepsilon$ \cite{fenichel1979geometric}, (\ref{GSPT_ALLEE_SYSTEMSLOW}) is transformed into an equivalent fast system, provided $\varepsilon \ne 0$,
\begin{equation} 
\label{GSPT_ALLEE_SYSTEMFAST}
    \left\{\begin{aligned}
  &\frac{\textrm{d}u}{\textrm{d}\xi}&&=w,\\
  &\frac{\textrm{d}w}{\textrm{d}\xi}&&=v+F(u),\\
  &\frac{\textrm{d}p}{\textrm{d}\xi}&&=\varepsilon R(u),\\
  &\frac{\textrm{d}v}{\textrm{d}\xi}&&=\varepsilon (p+cu).
\end{aligned}\right.
\end{equation}
The three fixed points\footnote{Even though the fixed points are independent of $\varepsilon$, we use the subscript $\varepsilon$ to indicate that these are fixed points of the full four-dimensional systems (\ref{GSPT_ALLEE_SYSTEMSLOW}) and (\ref{GSPT_ALLEE_SYSTEMFAST}).} of the two equivalent systems (\ref{GSPT_ALLEE_SYSTEMSLOW}) and (\ref{GSPT_ALLEE_SYSTEMFAST}) are
\begin{equation}
\label{FIX}
P^{0}_\varepsilon=(0,0,0,-F(0)),\,\, P^{1}_\varepsilon=(1,0,-c,-F(1)),\,\,
P^{A}_\varepsilon=(A,0,-cA,-F(A)),
\end{equation}
and we are interested in heteroclinic orbits connecting $P^{1}_\varepsilon$ with $P^{0}_\varepsilon$ as these correspond to travelling wave solutions of (\ref{RDE_biology3}) that asymptote to $1$ as $x\to-\infty$ and to $0$ as $x\to\infty$. 
Note that due to the symmetry
$
(w,p,z,c) \mapsto (-w,-p,-z,-c)
$
of system \eqref{GSPT_ALLEE_SYSTEMSLOW}, the existence of a heteroclinic orbit connecting $P^{1}_\varepsilon$ with $P^{0}_\varepsilon$ also implies the existence of a heteroclinic orbit connecting $P^{0}_\varepsilon$ with $P^{1}_\varepsilon$ and this latter orbit corresponds to a travelling wave solution of~(\ref{RDE_biology3}) that asymptotes to $0$ as $x\to-\infty$ and to $1$ as $x\to\infty$ and moves in the opposite direction.

The characteristic equation of the Jacobian of (\ref{GSPT_ALLEE_SYSTEMFAST}) is given by
\begin{equation}
    \label{eigenvalue1_1}
    \tau^4-D(u)\tau^2-\varepsilon\tau c-\varepsilon^2R'(u)=0,
\end{equation}
where we used that $F'(u)=D(u)$ and observe that $u=0,1$ or $A$ at a fixed point.
Upon substituting a regular expansion $\tau=\tau_0+\varepsilon\tau_1+\mathcal{O}(\varepsilon^2)$ into (\ref{eigenvalue1_1}),
we obtain an expansion for the four eigenvalues of the Jacobian 
\begin{equation}
    \label{eigs22}
    \begin{aligned}
    &\tau_{1}^{\pm}(u)=\frac{-c\pm\sqrt{c^2-4D(u)R'(u)}}{2D(u)}\varepsilon+\mathcal{O}(\varepsilon^2),\\ &\tau_{2}^{\pm}(u)=\pm\sqrt{D(u)}+\frac{c}{2D(u)}\varepsilon+\mathcal{O}(\varepsilon^2).
    \end{aligned}
\end{equation}
At $P^{0}_\varepsilon$, $R'(0)<0$, $D(0)>0$, thus $\tau_{1,2}^{+}(0)>0$, $\tau_{1,2}^{-}(0)<0$. 
Similarly, at $P^{1}_\varepsilon$, $R'(1)<0$, $D(1)>0$ and $\tau_{1,2}^{+}(1)>0$, $\tau_{1,2}^{-}(1)<0$.
That is, both the stable and unstable manifolds $P^{0,1}_\varepsilon$ are two-dimensional. At $P^{A}_\varepsilon$, the stable and unstable manifolds depend on the sign of $c$. If $c>0$, the stable manifold of $P^{A}_\varepsilon$ is three-dimensional and the unstable manifold of $P^{A}_\varepsilon$ is one-dimensional, while the situation for the stable and unstable manifolds of $P^{A}_\varepsilon$ is the opposite for $c<0$. For $c=0$, that is, for a standing wave, we again have that the stable and unstable manifold of $P^{A}_\varepsilon$ are two-dimensional. 

While the slow system (\ref{GSPT_ALLEE_SYSTEMSLOW}) and the fast system (\ref{GSPT_ALLEE_SYSTEMFAST}) are equivalent for $\varepsilon\ne0$, they have different singular limits. The singular limit of the fast system, that is, the layer problem, describes the dynamics near the shock and the fast variables $(u,w)$ will change significantly here while the slow variables $(p,v)$ are to leading order constant. In contrast, the singular limit of the slow system, that is, the reduced problem, describes the dynamics away from the shock and here the fast variables will be slaved to the slow variables. 

\subsection{Layer problem}
\begin{landscape}
\begin{figure}
\centering
\vspace{-1.5cm}
\includegraphics[]{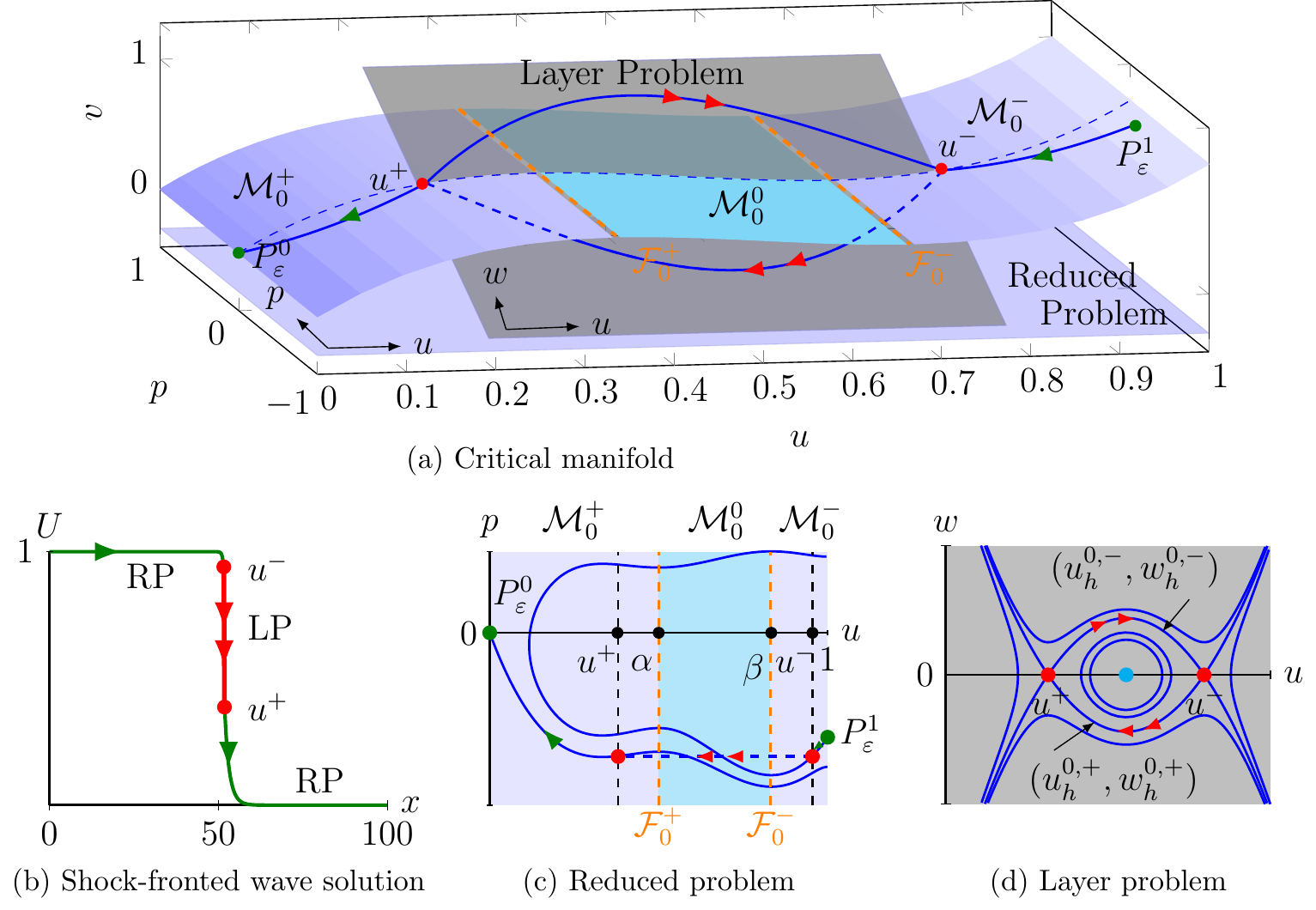}
\caption{(a) A projection of the four-dimensional phase plane of \eqref{GSPT_ALLEE_SYSTEMSLOW} and the critical manifold $\mathcal{M}_0$. A shock-fronted travelling wave solution of \eqref{RDE_biology3} as shown in (b) corresponds to a heteroclinic orbit (indicated in blue in (a)) that starts at $P^1_\varepsilon$ on the normally hyperbolic branch $\mathcal{M}_0^-$ of $\mathcal{M}_0$ and that follows the dynamics of the reduced problem (RP), whose projection on the $(u,p)$-plane is shown in (c), before it jumps to the other normally hyperbolic branch $\mathcal{M}_0^+$ according to the dynamics of the layer problem (LP). The projection of the layer dynamics on the $(u,w)$-plane, since $p$ and $v$ are constant, is shown in (d) and the two blue curves connecting $u^-$ and $u^+$ in (a) correspond to the heteroclinic orbits $\left(u^{0,\pm}_h,w^{0,\pm}_h\right)$ in (d). On $\mathcal{M}_0^+$, the heteroclinic orbit again follows the dynamics of the reduced problem and asymptotes to $P^0_\varepsilon$. A shock-fronted travelling wave solution is thus composed by orbits in the reduced problem and the layer problem as indicated in~(b).  } 
\label{PictureCriticalManifold}
\end{figure}
\end{landscape}

The layer problem is obtained by letting $\varepsilon\to0$ in the fast system (\ref{GSPT_ALLEE_SYSTEMFAST}), which gives
\begin{equation}
\label{GSPT_ALLEE_SYSTEMlayer_problem}
    \left\{\begin{aligned}
  &\frac{\textrm{d}u}{\textrm{d}\xi}&&=w,\\
  &\frac{\textrm{d}w}{\textrm{d}\xi}&&=v+F(u),
\end{aligned}\right.
\end{equation}
as well as $\textrm{d}p/\textrm{d}\xi=0$ and $\textrm{d}v/\textrm{d}\xi=0$, that is, $(p,v)\in\mathbb{R}^2$ are constants in (\ref{GSPT_ALLEE_SYSTEMlayer_problem}). The union of fixed points of (\ref{GSPT_ALLEE_SYSTEMlayer_problem})
\begin{equation}
    \label{CRIT}
    \mathcal{M}_0:=\{(u,w,p,v)\in\mathbb{R}^4:w=0, \ F(u)=-v\},
\end{equation} 
forms a two-dimensional invariant manifold, which is the so-called critical manifold \cite{jones1995geometric}, see Figure \ref{PictureCriticalManifold}.

The Jacobian of (\ref{GSPT_ALLEE_SYSTEMlayer_problem}) is
\begin{equation}
\nonumber
   J=\left(\begin{matrix}
   0 & 1 \\
   D(u) & 0
  \end{matrix}\right),
\end{equation}
with eigenvalues 
\begin{equation}
\label{eigenvalues1}
  \tau_{\pm}(u,w)=\pm\sqrt{D(u)}.
\end{equation}
Therefore, the manifold $\mathcal{M}_0$ loses normal hyperbolicity when $D(u)\leq0$, that is, for $u\in[\alpha,\beta]$ the eigenvalues $\tau_{\pm}$ (\ref{eigenvalues1}) are purely imaginary. As such, we split the critical manifold $\mathcal{M}_0$ into two two-dimensional normally hyperbolic saddle-type branches
\begin{equation}
    \nonumber
    \begin{aligned}
    &\mathcal{M}_0^{+}:=\{(u,w,p,v)\in\mathbb{R}^4:w=0, \ F(u)=-v, \ u\in[0,\alpha)\},\\
    &\mathcal{M}_0^{-}:=\{(u,w,p,v)\in\mathbb{R}^4:w=0, \ F(u)=-v, \ u\in(\beta,1]\},
    \end{aligned}
\end{equation} 

\noindent  
a two-dimensional not normally hyperbolic centre-type branche
\begin{equation}
    \nonumber
    \begin{aligned}
    &\mathcal{M}_0^{0}:=\{(u,w,p,v)\in\mathbb{R}^4:w=0, \ F(u)=-v, \ u\in(\alpha,\beta)\},\\
    \end{aligned}
\end{equation} 
and the two one-dimensional boundary sets
\begin{equation}
    \nonumber
    \begin{aligned}
    &\mathcal{F}_0^{+}:=\{(u,w,p,v)\in\mathbb{R}^4:w=0, \ F(u)=-v, \ u=\alpha\},\\
    &\mathcal{F}_0^{-}:=\{(u,w,p,v)\in\mathbb{R}^4:w=0, \ F(u)=-v, \ u=\beta\}.\\
    \end{aligned}
\end{equation} 

The layer problem \eqref{GSPT_ALLEE_SYSTEMlayer_problem} describes the dynamics near the shock away from the critical manifold.   It is a Hamiltonian system and, as such, we are looking for a heteroclinic orbit connecting $\mathcal{M}_0^{-}$ with $\mathcal{M}_0^{+}$. 
The Hamiltonian of (\ref{GSPT_ALLEE_SYSTEMlayer_problem}) is given by
\begin{equation*}
\label{hamiltonian_sec2}
    H(u,w)=-\frac{1}{2}w^2+G(u)+vu,
\end{equation*} 
where $G(u)= \int F(u) \,\textrm{d}u$. 
\begin{figure}
\centering
\includegraphics[width=\textwidth]{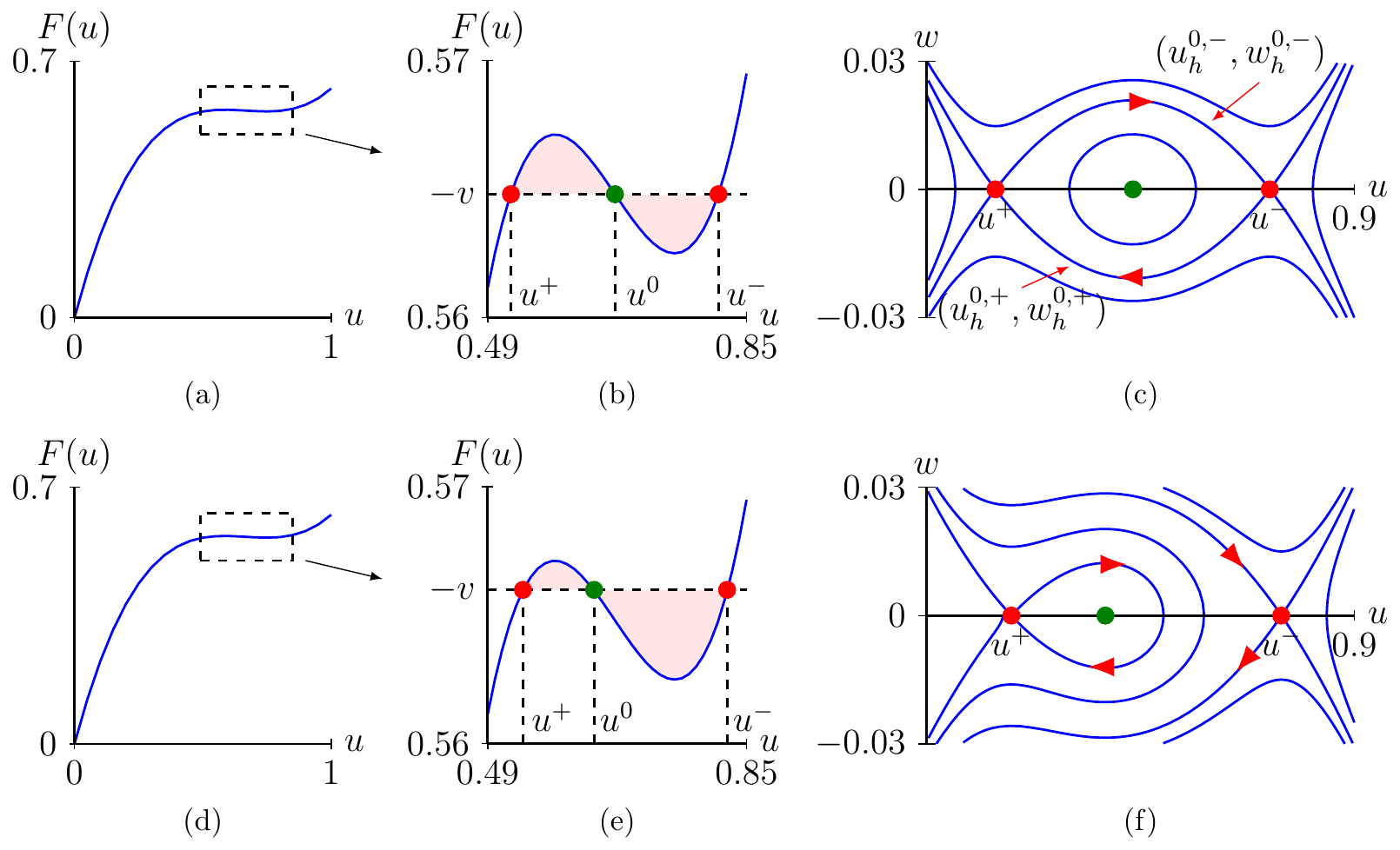}
\caption{(a) \& (b) The case when $u^+$ and $u^-$ satisfy the equal area rule (\ref{Equalarearule_ALLEE1}) with $D(u)=6(u-7/12)(u-3/4)$ leading to $F(u)=2u^3-4u^2+21u/8$ and $v=-61/108$. (c) The related phase plane of (\ref{GSPT_ALLEE_SYSTEMlayer_problem}) including the two heteroclinic orbits $\left(u_h^{0,\pm},w_h^{0,\pm}\right)$. (d) \& (e) The case when $u^+$ and $u^-$ do not satisfy the equal area rule with $v=-163/288$. (f) The related phase plane.} 
\label{Du&fu}
\end{figure}
Any solution is confined to a level set of the Hamiltonian and we have that 
\begin{equation}
\nonumber
    G(u^+)+vu^+=G(u^-)+vu^-,
\end{equation}
where $u^\pm \in \mathcal{M}_0^{\pm}$ are the end-points of the heteroclinic orbit such that $0< u^+<\alpha<\beta<u^-<1$.
This is equivalent to the integral equation
\begin{equation}
\label{Equalarearule_ALLEE1}
    \int_{u^+}^{u^-}\left(F(u)+v\right)\textrm{d}u=0,
\end{equation}
which is the well-known \emph{equal area rule}, see, for example, \cite{witelski1995shocks}. Recall that $F(u) = \int D(u) \,\textrm{d}u$ and $F(u)$ thus has an integration constant. Therefore, for a specific $F(u)$ the value of $v$ satisfying the equal area rule (\ref{Equalarearule_ALLEE1}) is unique. In~\ref{A:HET} we show that 
(\ref{GSPT_ALLEE_SYSTEMlayer_problem}) supports two heteroclinic orbits connecting $(u^+,0)$ and $(u^-,0)$ and
these heteroclinic orbits $\left(u_h^{0,\pm},w_h^{0,\pm}\right)$ are given by
\begin{equation*}
\begin{aligned}
    &u_{h}^{0,\pm}(\xi)=\frac{u^-+u^+}{2}\pm\frac{u^--u^+}{2}\tanh{\left(-\frac{a(u^--u^+)}{2}\xi\right)},\\
    &w_{h}^{0,\pm}(\xi)=\mp\frac{a(u^--u^+)^2}{4}\text{sech}^2{\left(-\frac{a(u^--u^+)}{2}\xi\right)},
\end{aligned}
\end{equation*}
where $a=\sqrt{(D_i-D_g)/2}$ and we recall that $u^->u^+$ by construction. See Figure \ref{Du&fu}. 
\subsection{Reduced problem}
The reduced problem is obtained from (\ref{GSPT_ALLEE_SYSTEMSLOW}) by letting $\varepsilon\to0$, which gives
\begin{equation}
\label{GSPT_ALLEE_SYSTEMreduced_problem}
    \left\{\begin{aligned}  
  &\frac{\textrm{d}p}{\textrm{d}z}&&=R(u),\\
  &\frac{\textrm{d}v}{\textrm{d}z}&&=p+cu,
\end{aligned}\right.
\end{equation}
and the two algebraic constraints $w=0$ and $v+F(u)=0$. Hence, (\ref{GSPT_ALLEE_SYSTEMreduced_problem}) simplifies to
\begin{equation}
\label{GSPT_ALLEE_SYSTEMreduced problem2_2}
    \left\{\begin{aligned}
  -D(u)&\frac{\textrm{d}u}{\textrm{d}z}&&=p+cu,\\
  &\frac{\textrm{d}p}{\textrm{d}z}&&=R(u).\\
\end{aligned}\right.
\end{equation}
Morever, since $w=0$ and $F(u)=-v$, \eqref{GSPT_ALLEE_SYSTEMreduced problem2_2} governs the flow on the critical manifold $\mathcal{M}_0$.  
The reduced problem is singular along the two lines $u=\alpha$ and $u=\beta$ since $D(\alpha)=D(\beta)=0$. Therefore, we transform (\ref{GSPT_ALLEE_SYSTEMreduced problem2_2}) into a desingularised system\footnote{Deriving the desingularised system from (\ref{GSPT_ALLEE_SYSTEMreduced problem2_2}) is, strictly speaking, not necessary for our analysis as we only need to consider $(\ref{GSPT_ALLEE_SYSTEMreduced problem2_2})$ on $\mathcal{M}^{\pm}_0$ away from $\alpha$ and $\beta$ since the heteroclinic orbit of the layer problem jumps from $\mathcal{M}^{-}_0 \ni u^- > \beta$ 
to $\in\mathcal{M}^+_0 \ni u^+ < \alpha$, see Figures \ref{PictureCriticalManifold} and \ref{Du&fu}. However, the desingularised system is more amenable to analysis and we thus study the dynamics of this desingularised system. 
} by using a stretched variable $\textrm{d}\psi=\textrm{d}z/D(u)$ \cite{Aronson1980161,2019arXiv190310090L}
\begin{equation}
\label{GSPT_ALLEE_SYSTEMreduced problem3_2}
    \left\{\begin{aligned}
  &\frac{\textrm{d}u}{\textrm{d}\psi}&&=-p-cu,\\
  &\frac{\textrm{d}p}{\textrm{d}\psi}&&=D(u)R(u).\\
\end{aligned}\right.
\end{equation}
It is important to note that, while the stretching changes the speed along a trajectory in a nonlinear fashion, the trajectories of the phase portraits of the reduced problem (\ref{GSPT_ALLEE_SYSTEMreduced problem2_2}) and the desingularised problem (\ref{GSPT_ALLEE_SYSTEMreduced problem3_2}) are the same. However, the orientation along a trajectory is reversed for $u\in(\alpha,\beta)$ as $D(u)<0$.

System (\ref{GSPT_ALLEE_SYSTEMreduced problem3_2}) has five fixed points $(0,0),\ (\alpha,-c\alpha),\ (A,-cA),\ (\beta,-c\beta)$ and $(1,-c)$. 
The eigenvalues  and eigenvectors of the Jacobian of (\ref{GSPT_ALLEE_SYSTEMreduced problem3_2}) are given by
\begin{equation}
    \nonumber
    \tau_{\pm}=\frac{-c\pm\sqrt{c^2-4(D(u)R'(u)+D'(u)R(u))}}{2},\quad E_{\pm}=(1,-\tau_{\mp}).
\end{equation} 
If we let $\chi_1$ be the minimum of the set $\{A,\alpha,\beta\}$, $\chi_3$ its maximum, and $\chi_2$ the remaining element, then the characteristics of $D(u)$ and $R(u)$ yield the following results:
\begin{itemize}
\item $(0,0), (1,-c)$ and $(\chi_2,-c \chi_2)$ are saddles; and
\item 
$(\chi_i,-c\chi_i)$, $i \in \{1,3\}$, is a(n) 
\begin{itemize}
\item
stable node for $c>2\sqrt{D(\chi_i)R'(\chi_i)}$;
\item
stable spiral for $0<c<2\sqrt{D(\chi_i)R'(\chi_i)}$; 
\item centre for $c=0$;
\item unstable spiral for $-2\sqrt{D(\chi_i)R'(\chi_i)}<c<0$; and
\item unstable node for $c<-2\sqrt{D(\chi_i)R'(\chi_i)}$.
\end{itemize}
\end{itemize}

\subsection{The construction of the heteroclinic orbit in the singular limit}
Since the fixed points $P^{0,1}_\varepsilon$ \eqref{FIX} are on the normally hyperbolic branches $\mathcal{M}_0^\pm$ of the critical manifold $\mathcal{M}_0$, a shock-fronted travelling wave solution to~\eqref{RDE_biology3} corresponds to a heteroclinic orbit of \eqref{GSPT_ALLEE_SYSTEMSLOW} that, to leading order, starts on $\mathcal{M}_0^-$, follows the dynamics of the reduced problem \eqref{GSPT_ALLEE_SYSTEMreduced problem2_2} before it jumps, according to the layer dynamics \eqref{GSPT_ALLEE_SYSTEMlayer_problem}, to the other normally hyperbolic branch $\mathcal{M}_0^+$ on which it asymptotes to $P^{0}_\varepsilon$ following the dynamics of \eqref{GSPT_ALLEE_SYSTEMreduced problem2_2} again. In particular, if we split the spatial domain $z \in (-\infty,\infty)$ into three parts: 
\begin{align}
\label{split}
\begin{aligned}
z \in I_s^- &:= (-\infty,-\sqrt{\varepsilon}+z^*)\,, \quad
z \in I_f := [-\sqrt{\varepsilon}+z^*, \sqrt{\varepsilon}+z^*]\,, \\
z \in I_s^+ &:= (\sqrt{\varepsilon}+z^*, \infty)\,,
\end{aligned}
\end{align} 
then the heteroclinic orbit is, to leading order, on $\mathcal{M}_0^\pm$ and governed by the reduced problem \eqref{GSPT_ALLEE_SYSTEMreduced problem2_2} for $z \in  I_s^\pm$, while it is, to leading order, governed by the layer problem  \eqref{GSPT_ALLEE_SYSTEMlayer_problem} for $z \in  I_f$, see Figure~\ref{PictureCriticalManifold}. Note that due to translation invariance of (\ref{GSPT_ALLEE_ODE1_partial}) we can, without loss of generality, set $z^*=0$ in \eqref{split}.

\begin{figure}
\vspace{-2.5cm}
\centering
\includegraphics[]{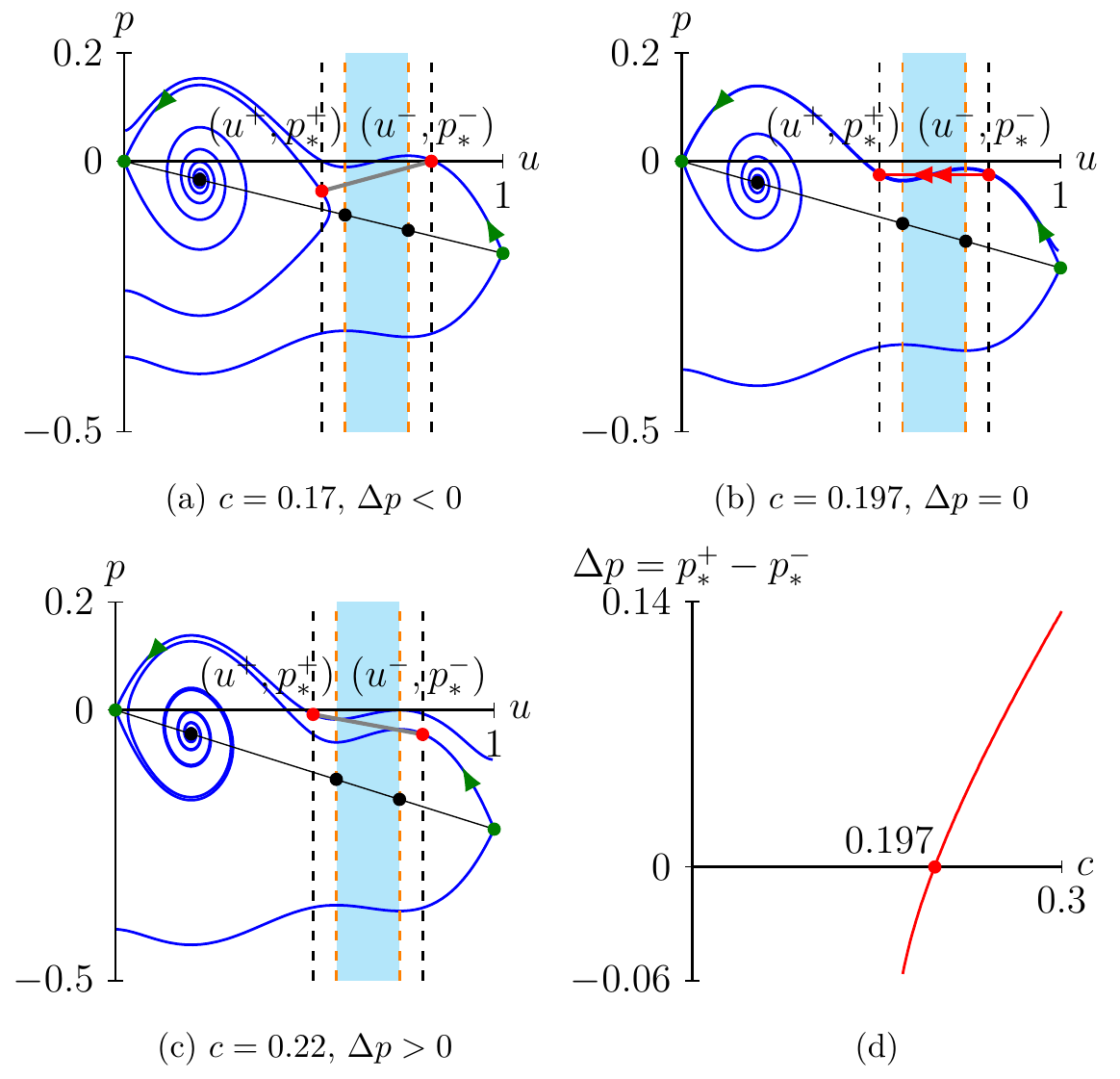}\\
\caption{(a) -- (c) Phase planes of the desingularised system (\ref{GSPT_ALLEE_SYSTEMreduced problem3_2}) for different values of $c$ with $D(u)=6(u-7/12)(u-3/4)$ and $R(u)=5u(1-u)(u-1/5)$. The green and black dots are fixed points and red points are the endpoints $u^\pm$ of the fast jump (as derived from the layer problem). The black solid lines are the nullclines $p=-cu$. The red straight solid line in (b) indicates the shock from $(u^-,p_*^{-})$ to $(u^+,p^{+}_*)$, while the grey straight lines in (a) \& (c) cannot 
lead to shocks as $p_*^{+}\ne p_*^{-}$. With $c = 0.197$, it leads to a feasible desingularised system where $p_*^{+}=p_*^{-}$. (d) The difference $\Delta p = p_*^+-p_*^-$ as a function of $c$ shows that $\Delta p=0$ for $c=c_0\approx0.197$.}
\label{Panel_Martin_relation1}
\end{figure}

Since $w=0$ and $F(u)=-v $ on the critical manifold, the fixed points $(0,0)$ and $(1,-c)$ of the reduced problem \eqref{GSPT_ALLEE_SYSTEMreduced problem2_2} correspond to $P^{0}_\varepsilon$ and $P^{1}_\varepsilon$, respectively. Furthermore, the analysis of the layer problem (\ref{GSPT_ALLEE_SYSTEMlayer_problem}) -- which is independent of the speed $c$ -- indicates there may exist shocks with endpoints $u^-$ ($>\beta$) and $u^+$ ($<\alpha$). Consequently, if there exists a shock-fronted travelling wave solution of~(\ref{RDE_biology3}) with a shock from $u^-$ to $u^+$, it relates to two trajectories in 
system~(\ref{GSPT_ALLEE_SYSTEMreduced problem2_2}), see also Figure~\ref{PictureCriticalManifold}. These, in turn, relate to two corresponding trajectories in the desingularised system (\ref{GSPT_ALLEE_SYSTEMreduced problem3_2}). One is the unique trajectory $
\gamma^+$, for a given speed $c$, that starts on the line $\{(u^+,p^+), p^+ \in \mathbb{R}\}$ and approaches $(0,0)$ as $\psi\to \infty$, while the other one is the unique trajectory $
\gamma^-$ that arrives at the line $\{(u^-,p^-), p^- \in \mathbb{R}\}$ and approaches $(1,-c)$ as $\psi\to-\infty$.
Note that these unique trajectories can intersect the lines $\{(u^{\pm},p^\pm), p^\pm \in \mathbb{R}\}$ multiple times, see, for instance, Figure~\ref{Panel_Martin_relation1}a.  However, only the first intersections may lead to monotone travelling wave solutions. 
\noindent 
We are mainly interested in monotone travelling wave solutions since nonmonotonic travelling wave solutions are often PDE unstable \cite{VOLPERT}. Therefore, we only look for these first intersections.
As $p$ is a slow variable, it should, to leading order, hold constant at the endpoints of the shock ($\textrm{d}p/\textrm{d}\xi=0$ in the singular limit). Hence, we are interested in the speeds $c_0$ for which
the $p$-value of the trajectory $\gamma^-$ at $u^-$, say $p_*^-$, is the same as the $p$-value of the trajectory $\gamma^+$ at $u^+$, say $p_*^+$, see Figure~\ref{Panel_Martin_relation1}b.
These $c$-values  
determine the actual speed of the shock-fronted travelling wave solution.

As the stable and unstable manifolds of $(0,0)$ and $(1,-c)$ are algebraically too complicated to study analytically, we use numerical tools to detect the speeds leading to a feasible desingularised system (\ref{GSPT_ALLEE_SYSTEMreduced problem3_2}). In particular, we use the function \textit{ode45} in MATLAB to obtain the phase plane of (\ref{GSPT_ALLEE_SYSTEMreduced problem3_2}) and then calculate $\Delta p :=p_*^{+}-p_*^{-}$ for different speeds $c$. Note that we locate the initial points of trajectories approaching $(0,0)$ or $(1,-c)$ with a small step along their eigenvectors. We find the crossing point of the trajectory leaving from $(1,-c)$ and the straight line $u=u^-$ as $(u^-,p_*^{-})$ and the crossing point of the trajectory arriving at $(0,0)$ and the straight line $u=u^+$ as $(u^+,p_*^{+})$. Finally, we calculate $\Delta p$ as function of $c$.

As shown in Figure \ref{Panel_Martin_relation1}, for a given prototypical $D(u)=6(u-7/12)(u-3/4)$ and $R(u)=5u(1-u)(u-1/5)$, the difference between two $p$-values at $u^{\pm}$ is zero when $c=c_0\approx0.197$, that is, the phase plane of (\ref{GSPT_ALLEE_SYSTEMreduced problem3_2}) aligns the endpoints of the shock $u^-\to u^+$ when $c=c_0\approx0.197$. Thus, in the singular limit $\varepsilon \to 0$, we expect that \eqref{RDE_biology3}, with the given $D(u)$ and $R(u)$, supports a shock-fronted right-travelling wave solution with speed $c=c_0\approx0.197$, see Figure \ref{Panel_TW1}b.

\begin{figure}
\centering
\includegraphics[width=\textwidth]{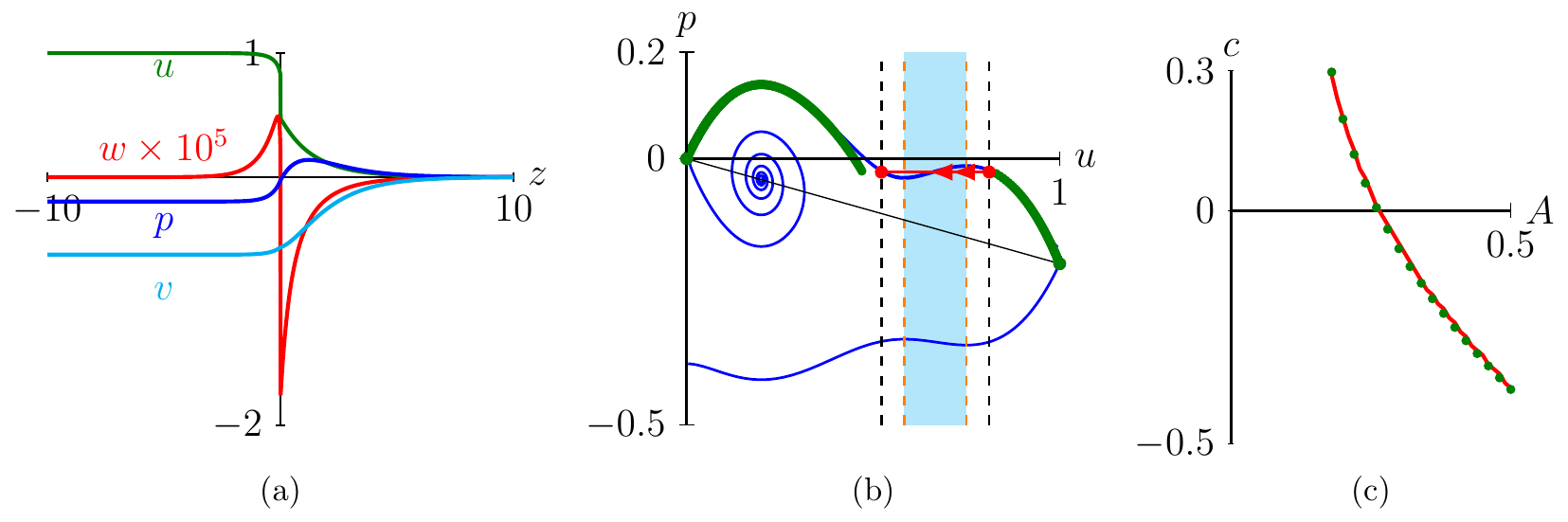}
\caption{The comparison between numerical and analytical results. (a) Numerical simulation of \eqref{GSPT_ALLEE_SYSTEMSLOW} with $D(u)=6(u-7/12)(u-3/4)$, $R(u)=5u(1-u)(u-1/5)$ and $\varepsilon=10^{-5}$ obtained using Matlab ODE solver `bvp4c'. Boundary conditions are $u(-L)=1, w(-L)=0, w(L)=0, p(-L)=-c, v(L)=-F(0)$ and $q(-L)=-5$ (where $\textrm{d}q/\textrm{d}z=u$). The numerically obtained speed is $c=0.196$, which is close to the analytically obtained speed $c_0=0.197$. 
(b) The numerically obtained $(u,p)$-trajectory in terms of the stretched variable $\psi$ (green line) superimposed onto the phase plane of the desingularised system (\ref{GSPT_ALLEE_SYSTEMreduced problem3_2}) with $c=c_0=0.197$. 
(c) The numerically (green) and analytically (red) obtained speeds for a varying reaction term of the form $R(U)=5U(1-u)(u-A)$ with $A \in (0, 0.5)$.
}
\label{Relation_Sec4_1}
\end{figure}

Due to the complexity of numerically simulating a singularly perturbed fourth-order PDE like \eqref{RDE_biology3}, we simulate solutions of the perturbed ODE system \eqref{GSPT_ALLEE_SYSTEMSLOW} with Matlab's ODE solver `bvp4c' and compare it with our analytical results from the singular limit.
With the diffusivity function and reaction term as above, the numerical results and analytical results coincide (to leading order), see Figure \ref{Relation_Sec4_1}. Furthermore, in Figure \ref{Relation_Sec4_1}c we compare the numerical and analytical speeds for reaction terms of the form $R(u)=5u(1-u)(u-A)$ with varying $A$. Again, the numerical and analytical speeds coincide (to leading order). 

\subsection{Persistence analysis}
For $c=c_0$, the orbit in the layer problem connecting $u^-$ to $u^+$ and the orbits in the reduced problem and desingularised problem connecting $1$ to $u^-$ and connecting $u^+$ to $0$ form a complete heteroclinic orbit connecting $1$~to~$0$ in the singular limit $\varepsilon\to0$. Below we will argue that such solution persists in the four-dimensional system (\ref{GSPT_ALLEE_SYSTEMSLOW}) for sufficiently small $
\varepsilon$, i.e. $0<\varepsilon\ll1$. 
Note that we do not present the full proof for the persistence claim -- which follows from geometric singular perturbation theory (GSPT) based on Fenichel's persistence theorems \cite{fenichel1979geometric, hek2010geometric,jones1995geometric} since $\mathcal{M}_0^{\pm}$ are normally hyperbolic -- because this is 
rather standard, but quite technical, at this stage. Instead, we provide some heuristic arguments for the persistence.

The endpoints of the heteroclinic orbit in the full system (\ref{GSPT_ALLEE_SYSTEMSLOW}) are $P_\varepsilon^0$ and $P_\varepsilon^1$ \eqref{FIX} and the heteroclinic orbit lies in the intersection of the two-dimensional stable manifold of $P_\varepsilon^{0}$, $\mathcal{W}^s(P_\varepsilon^{0})$, and the two-dimensional unstable manifold of $P_\varepsilon^{1}$, $\mathcal{W}^u(P_\varepsilon^{1})$, see \eqref{eigs22}. This intersection will generically not be transversal since the full system is four-dimensional, i.e. $2+2-1<4$. Therefore, we extend the full system (\ref{GSPT_ALLEE_SYSTEMSLOW}) to a five-dimensional system by appending it with an equation for the unknown speed 
$\{c'=0\}$.
That is, we threat $c$ as a variable and not as an unknown parameter.
In the extended system the heteroclinic orbit now lies in the intersection of the three-dimensional centre stable manifold $\mathcal{W}^{cs}(P_\varepsilon^{0})$ and the three-dimensional centre unstable manifold $\mathcal{W}^u(P_\varepsilon^{1})$ and this intersection will generically be transversal since the full system is five-dimensional, i.e. $3+3-1=5$. Typically, transversality follows from a Melnikov-type analysis \cite{hek2010geometric,Robinson1983,szmolyan1991transversal}. We decided to omit this calculation, but its proof is numerically verified in Figure~\ref{Panel_Martin_relation1}(b) and (d). As a result, and for sufficiently small $\varepsilon$, the heteroclinic orbit will persist with a nearby speed $c(\varepsilon)$, with $c(0)=c_0$, the speed found in the singular limit. 
Finally, recall that such a heteroclinic orbit corresponds to a shock-fronted travelling wave solution of \eqref{RDE_biology3}.

\section{Viscous relaxation}
\label{sec3}
In this section, we study shock-fronted travelling wave solutions in (\ref{Another1}) and we use similar mathematical techniques as in \S\ref{sec2} to obtain a three-dimensional singular perturbed slow-fast system. 
 The reduced problem is the same as in \S\ref{sec2}, however, it has different algebraic constraints. In contrast, the layer problem is different and only one-dimensional which leads to shocks with different characteristics. Since the methodology of the analysis is similar, we only present a succinct and brief derivation of the main results.\subsection{Preliminary observations}
The travelling wave solution of (\ref{Another1}) of interest here is a solution of 
\begin{equation}
\label{sec3_ODE1}
    \begin{aligned}
    \frac{\textrm{d}}{\textrm{d}z}\left(\varepsilon c\frac{\textrm{d}^2u}{\textrm{d}z^2}-cu-\frac{\textrm{d}}{\textrm{d}z}(F(u))\right)=R(u),
    \end{aligned}
\end{equation}
that asymptotes to one as $z \to -\infty$ and to zero as $z \to \infty$. Here,  $z:=x-ct$ is again the travelling wave coordinate. Next, with some abuse of notation, we define
\begin{equation}
\label{definationtoODE2}
    p:=\varepsilon c\frac{\textrm{d}^2u}{\textrm{d}z^2}-cu-\frac{\textrm{d}}{\textrm{d}z}(F(u)),\quad
    v:=\varepsilon  c\frac{\textrm{d}u}{\textrm{d}z}-F(u),
\end{equation}
and transform (\ref{sec3_ODE1}) into a three-dimensional singular perturbed slow-fast dynamical system
\begin{equation}
\label{Another4_1}
   \left\{\begin{aligned}
   \varepsilon&\frac{\textrm{d}u}{\textrm{d}z}&&=\frac{1}{c}(v+F(u)),\\
   &\frac{\textrm{d}p}{\textrm{d}z}&&=R(u),\\
   &\frac{\textrm{d}v}{\textrm{d}z}&&=p+cu,
\end{aligned}\right.
\end{equation} 
where $u\in\mathbb{R}$ is fast variable and $(p,v)\in\mathbb{R}^2$ are slow variables. By using a stretched variable $\xi=z/\varepsilon$, (\ref{Another4_1}) is transformed into an equivalent fast system, provided $\varepsilon\ne0$,
\begin{equation}
\label{Another4}
    \left\{\begin{aligned}
   &\frac{\textrm{d}u}{\textrm{d}\xi}&&=\frac{1}{c}(v+F(u)),\\
   &\frac{\textrm{d}p}{\textrm{d}\xi}&&=\varepsilon R(u),\\
   &\frac{\textrm{d}v}{\textrm{d}\xi}&&=\varepsilon (p+cu).
\end{aligned}\right.
\end{equation}
The fixed points of the two equivalent systems (\ref{Another4_1}) and (\ref{Another4}) are
\begin{equation}
    \nonumber
    Q_\varepsilon^{0}=(0,0,-F(0)),\quad Q_\varepsilon^{A}=(A,-cA,-F(A)),\quad
    Q_\varepsilon^{1}=(1,-c,-F(1)),
\end{equation}
and we are interested in heteroclinic orbits connecting $Q_\varepsilon^{0}$ with $Q_\varepsilon^{1}$.
The Jacobian of (\ref{Another4}) has three eigenvalues with the expansion of $\varepsilon$
\begin{equation}
    \nonumber
    \begin{aligned}
    &\tau_{1}^{\pm}(u)=\frac{-c\pm\sqrt{c^2-4D(u)R'(u)}}{2D(u)}\varepsilon+\mathcal{O}(\varepsilon^2),
    \\
    &\tau_{2}(u)=\frac{D(u)}{c}+\frac{c}{D(u)}\varepsilon+\mathcal{O}(\varepsilon^2).
    \end{aligned}
\end{equation}
At $Q^{0}_\varepsilon$, $R'(0)<0$, $D(0)>0$, thus, $\tau_{1}^{+}(0)>0$, $\tau_{1}^{-}(0)<0$ and $\tau_{2}(0)>0$. 
Similarly, at $Q^{1}_\varepsilon$, $R'(1)<0$, $D(1)>0$, thus, $\tau_{1}^{+}(1)>0$, $\tau_{1}^{-}(1)<0$ and $\tau_{2}(1)>0$.
That is, the stable manifolds of $Q^{0,1}_\varepsilon$ are one-dimensional and the unstable manifolds of $Q^{0,1}_\varepsilon$ are two-dimensional.  
At $Q^{A}_\varepsilon$, for positive speeds, the stable manifold is two-dimensional and the unstable manifold is one-dimensional; for negative speeds, the stable manifold is one-dimensional and the unstable manifold is two-dimensional.

\subsection{Layer problem}
\begin{figure}
\centering
\includegraphics[width=\textwidth]{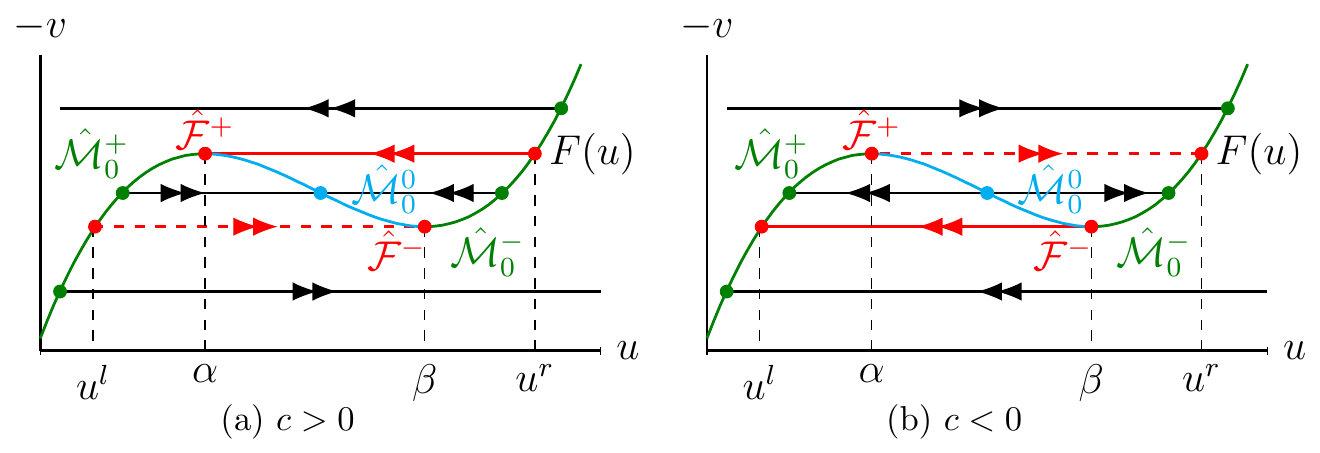}
\caption{(a) Connections from $\hat{\mathcal{M}}^{\pm}_0$ to $\hat{\mathcal{F}}_{\mp}$ for $c>0$. (b) These connections are reversed for $c<0$.
}
\label{Du&fu_more1}
\end{figure}

Letting $\varepsilon\to0$ in (\ref{Another4}) gives the layer problem
\begin{equation}
\label{Another_layer_problem}
    \frac{\textrm{d}u}{\textrm{d}\xi}=\frac{1}{c}(v+F(u)),
\end{equation}
and $\textrm{d}p/\textrm{d}\xi=0$ and ${\textrm{d}v}/{\textrm{d}\xi}=0$. Thus, we have a two-dimensional critical manifold
\begin{equation}
    \nonumber
    \hat{\mathcal{M}}_{0}:=\{(u,p,v)\in\mathbb{R}^3:F(u)=-v\}.
\end{equation}
Upon recalling that $F'(u)=D(u)$, we observe that the critical manifold loses normal hyperbolicity along the one-dimensional set
\begin{equation}
    \nonumber
    \hat{\mathcal{F}}:=\{(u,p,v)\in\hat{\mathcal{M}}_0:D(u)=0\},
\end{equation}
which has two branches
\begin{equation}
    \nonumber
    \hat{\mathcal{F}}=\hat{\mathcal{F}}^+\cup\hat{\mathcal{F}}^-:=\{(u,p,v)\in\hat{\mathcal{M}}_0:u=\alpha\}\cup\{(u,p,v)\in\hat{\mathcal{M}}_0:u=\beta\}.
\end{equation}
Thus, we split the critical manifold into five branches
$
    \hat{\mathcal{M}}_{0}=\hat{\mathcal{M}}_{0}^-\cup \hat{\mathcal{F}}^-\cup\hat{\mathcal{M}}_{0}^0\cup \hat{\mathcal{F}}^+\cup\hat{\mathcal{M}}_{0}^{+},$
with
\begin{equation}
    \nonumber
    \hat{\mathcal{M}}_{0}^+:=\{(u,p,v)\in\mathbb{R}^3: u<\alpha\},\quad
    \hat{\mathcal{M}}_{0}^-:=\{(u,p,v)\in\mathbb{R}^3: u>\beta\},
\end{equation}
repelling manifolds for $c>0$ and attracting manifolds for $c<0$. Similarly
\begin{equation}
    \nonumber
    \begin{aligned}
    &\hat{\mathcal{M}}_{0}^0=\{(u,p,v)\in\mathbb{R}^3: \alpha<u<\beta\},
    \end{aligned}
\end{equation}
is an attracting manifold for $c>0$ and an unstable manifold for $c<0$, see Figure \ref{Du&fu_more1}. 

Considering the stability of the different branches of critical manifold, there may exist connections between $\hat{\mathcal{M}}^{\pm}_0$ and $\hat{\mathcal{M}}^{0}_0$ and between $\hat{\mathcal{M}}^{\pm}_0$ and $\hat{\mathcal{F}}^{\mp}$.
In contrast to the previous section, we are now interested in connections between $\hat{\mathcal{M}}_0^{\pm}$ and $\hat{\mathcal{F}}^{\mp}$ since we are looking for travelling wave solutions that connect $u=0$ and $u=1$, and both of these points are on $\hat{\mathcal{M}}^{\pm}_0$. 
There are two ways to establish these connections. If $c>0$,
$\hat{\mathcal{M}}^{\pm}_0$ are repelling and $F(u)=-v$ has two non-repeating real roots $\beta$ and $u^l$ ($<\alpha$), or $u=\alpha$ and $u^r$ ($>\beta$). In this case, the related shocks are $u^r\to\alpha$ and $u^l\to\beta$, see Figure~\ref{Du&fu_more1}a. If $c<0$, $\hat{\mathcal{M}}^{\pm}_0$ are attracting. and the related shocks are in the opposite direction $\alpha\to u^r$ and $\beta\to u^l$, see Figure \ref{Du&fu_more1}b.

\subsection{Reduced problem} 
The reduced problem of (\ref{Another4_1}), obtained by letting $\varepsilon \to 0$, is the same as the reduced problem $(\ref{GSPT_ALLEE_SYSTEMreduced problem2_2})$ of the previous section and is given by
\begin{equation*}
    \left\{\begin{aligned}
  -D(u)&\frac{\textrm{d}u}{\textrm{d}z}&&=p+cu,\\
  &\frac{\textrm{d}p}{\textrm{d}z}&&=R(u).\\
\end{aligned}\right.
\end{equation*}
Similarly, its desingularised system\footnote{The desingularised system is required this time because we need to study the dynamics around $u=\alpha$ and $u=\beta$ where the reduced problem is singular.} is the same and given by
\begin{equation}
\label{GSPT_ALLEE_SYSTEMreduced problem3_3}
    \left\{\begin{aligned}
  &\frac{\textrm{d}u}{\textrm{d}\psi}&&=-p-cu,\\
  &\frac{\textrm{d}p}{\textrm{d}\psi}&&=D(u)R(u).\\
\end{aligned}\right.
\end{equation}
However, note that the slow variable $p$ is defined differently, see (\ref{definationtoODE1}) and (\ref{definationtoODE2}), and thus has a different meaning. 

\subsection{The construction of the heteroclinic orbit in the singular limit}
From the analysis of the layer problem~(\ref{Another_layer_problem}), the shocks $u^r\to\alpha$ and $u^l\to\beta$ have positive speeds, while the shocks in the opposite directions, $\alpha \to u^r$ and $\beta \to u^l$, have negative speed. 
The shocks $u^r\to\alpha$ and $\beta \to u^l$ potentially relate, in the singular limit, to trajectories of (\ref{Another4_1}) leaving from $u=1$ and arriving at $u=0$, that is, they have the asymptotic conditions $\lim_{z\to-\infty}u=1$ and $\lim_{z\to\infty}u=0$ we are interested in. In contrast, the shocks $u^l\to\beta$  and $\alpha \to u^r$ correspond to trajectories with the opposite asymptotic conditions $\lim_{z\to-\infty}u=0$ and $\lim_{z\to\infty}u=1$. 
Thus, we are interested in positive speeds $c$ for which there exist trajectories of the desingularised system $(\ref{GSPT_ALLEE_SYSTEMreduced problem3_3})$ that connect 
$(1,-c)$ with $(u^r, p_*)$ and
$(\alpha, p_*)$ with $(0,0)$ (both in forward $\psi$). Similarly, we are interested in negative speeds $c$ for which there exist trajectories of the desingularised system $(\ref{GSPT_ALLEE_SYSTEMreduced problem3_3})$ that connect $(1,-c)$ with $(\beta, p_*)$ and $(u^l, p_*)$ with $(0,0)$.

\begin{figure}
\centering
\includegraphics[width=\textwidth]{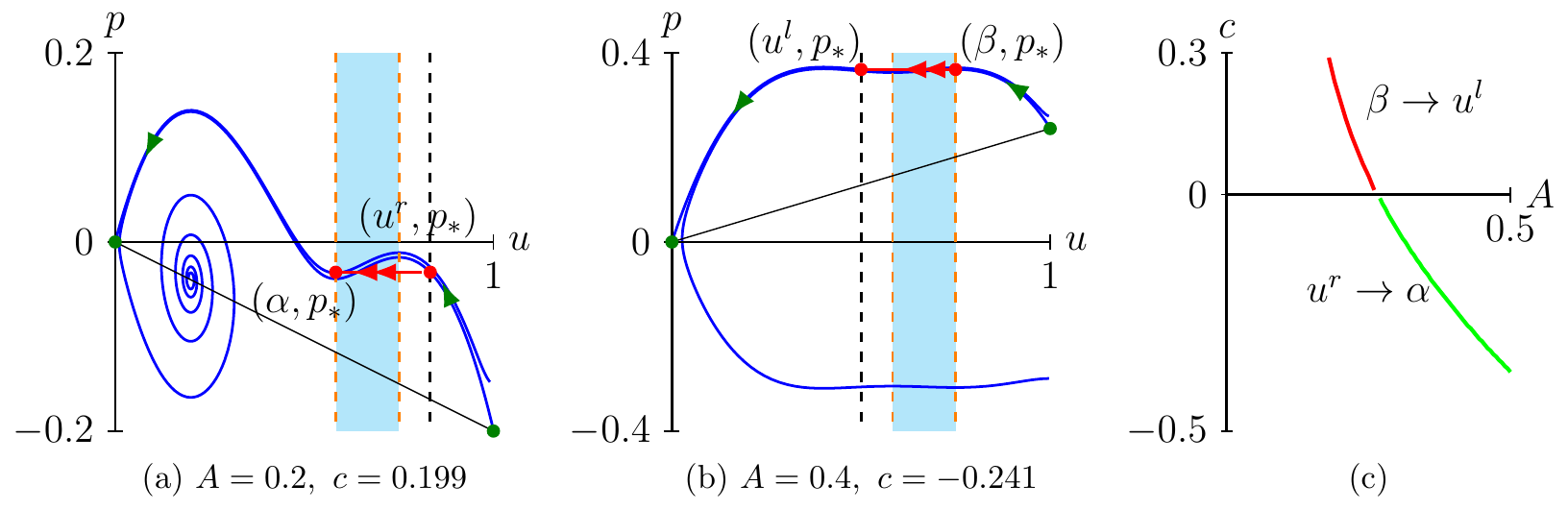}
\caption{(a) and (b) Phase planes of the desingularised system (\ref{GSPT_ALLEE_SYSTEMreduced problem3_2}) with $D(u)=6(u-7/12)(u-3/4)$ and $R(u)=5u(1-u)(u-A)$ with $A=0.2$ (a) and $A=0.4$ (b). In the former case, we observe a shock $u^r\to\alpha$ with a positive speed $c=c_0=0.199$, while in the latter case we have a shock $\beta\to u^l$ with a negative speed $c=c_0=-0.241$. (c) The change of speed as function of  $A$ where the line with positive speed represents shocks $u^r\to\alpha$ and the line with negative speed represents shocks $\beta\to u^l$. We remark that we could not find an Allee type $R(u)$ for which both types of travelling wave solutions exist simultaneously. } 
\label{Phaseplane1_2}
\end{figure}

Following the same procedure as in the previous section using {\emph{ode45}} in MATLAB, we can now construct orbits of the correspond to heteroclinic orbits in the singular limit of (\ref{Another4_1}), and thus to shock-fronted travelling wave solutions of \eqref{Another1}. See Figure~\ref{Phaseplane1_2} for two prototypical examples of these orbits. One corresponding to a  shock-fronted travelling wave solution with positive speed and one with negative speed.

\subsection{Persistence analysis}
To show that these singular orbits indeed persist for $\varepsilon$ sufficient small, and thus correspond to shock-fronted travelling wave solutions of \eqref{Another1}, we have to proceed in a similar fashion as in the previous section and extend the full three-dimensional system \eqref{Another4_1} with an equation for the speed $\{c'=0\}$ (since the stable and unstable manifold $\mathcal{W}^{s,u}(Q_\varepsilon^{0,1})$ are respectively one and two-dimensional and $1+2-1<3$) such that transversality is generically possible. 
Transversality again follows from a Melnikov-type argument, but  we have to 
extend Fenichel theory near the regular fold point $\hat{\mathcal{F}}^{\pm}$, where the critical manifold loses normal hyperbolicity -- one of the necessary conditions for Fenichel's persistence theorems.
This way, we can show that the orbits persist even though in the singular limit we leave, or arrive at, the  
critical manifold at a fold point $\hat{\mathcal{F}}^{\pm}$.
We decided to not go into the details of this analysis and refer to \cite{BECK}, and references therein, instead, for an outline how the persistence of these singular orbits can be shown. In the end, this shows the persistence of the heteroclinic orbit for sufficiently small $\varepsilon$ and with nearby speed $c(\varepsilon)$, with $c(0)=c_0$, the speed found in the singular limit.

\section{Summary, discussion and outlook}
\label{summary}
In this article, we studied shock-fronted travelling wave solutions supported by the RDE (\ref{RDE_biology1}) with a convex nonlinear diffusivity function $D(U)$ (\ref{D(u)}) that is negative for $U \in (\alpha, \beta)$ \eqref{AB}, and with an Allee-type reaction-term $R(U)$ (\ref{AlleeR(u)}). This RDE with forward-backward diffusion was previously derived by \cite{johnston2017co} from a lattice-based stochastic model modelling a population of individuals and groups that can undergo movement, birth and death events to describe the its macroscopic behaviour. 
We studied the RDE by adding two different small regularisations; a non-local regularisation $-\varepsilon^2 \partial^4 U/\partial x^4$, with $\varepsilon$ small, see \eqref{RDE_biology3} and \S\ref{sec2}, and a viscous relaxation $\varepsilon \partial^3 U/(\partial x^2 \partial t)$, see \eqref{Another1} and \S\ref{sec3}. Note that in the singular limit $
\varepsilon \to 0$ both PDEs reduce to (\ref{RDE_biology1}). 

These two regularisations allowed us to use a dynamical systems approach to study the shock-fronted travelling wave solutions. In particular, for the non-local regularisation the PDE \eqref{RDE_biology3} could be reduced to a singularly perturbed four-dimensional system of ODEs (\ref{GSPT_ALLEE_SYSTEMSLOW}). 
As the regularisation term is assumed to be small there is a scale separation in this system of ODEs. This allowed for a further reduction by investigating (\ref{GSPT_ALLEE_SYSTEMSLOW}) singular limit in the fast and slow scaling. The singular limit in the fast scaling, called the layer problem~\eqref{GSPT_ALLEE_SYSTEMlayer_problem}, described the dynamics near the shock of a shock-fronted travelling wave solutions, and was a two-dimensional Hamiltonian system independent of the speed $c$, see Figure~\ref{Du&fu}. 
The singular limit in the slow scaling, called the reduced problem \eqref{GSPT_ALLEE_SYSTEMreduced_problem}, was a singular two-dimensional system of ODEs. It is constraint to the critical manifold $
\mathcal{M}_0$ \eqref{CRIT} and described the dynamics away from the shock. Note that we use MATLAB to investigate the reduced problem as it is algebraically too involved to determine the sought after trajectories. A shock-fronted travelling wave solution can now be constructed, in the singular limit, upon concatenating the three parts of the solution, see Figures~
\ref{PictureCriticalManifold} and~\ref{Panel_Martin_relation1}. Subsequently, GSPT can be used to show that the solution persists for sufficiently small $\varepsilon$. Note that the details of this final calculation were omitted, instead it was shown that the dynamics of the full ODE (\ref{GSPT_ALLEE_SYSTEMSLOW}) agrees with the obtained results in the singular limit, see Figure~\ref{Relation_Sec4_1}.

\begin{figure}
\vspace{-2cm}
\centering
\includegraphics[]{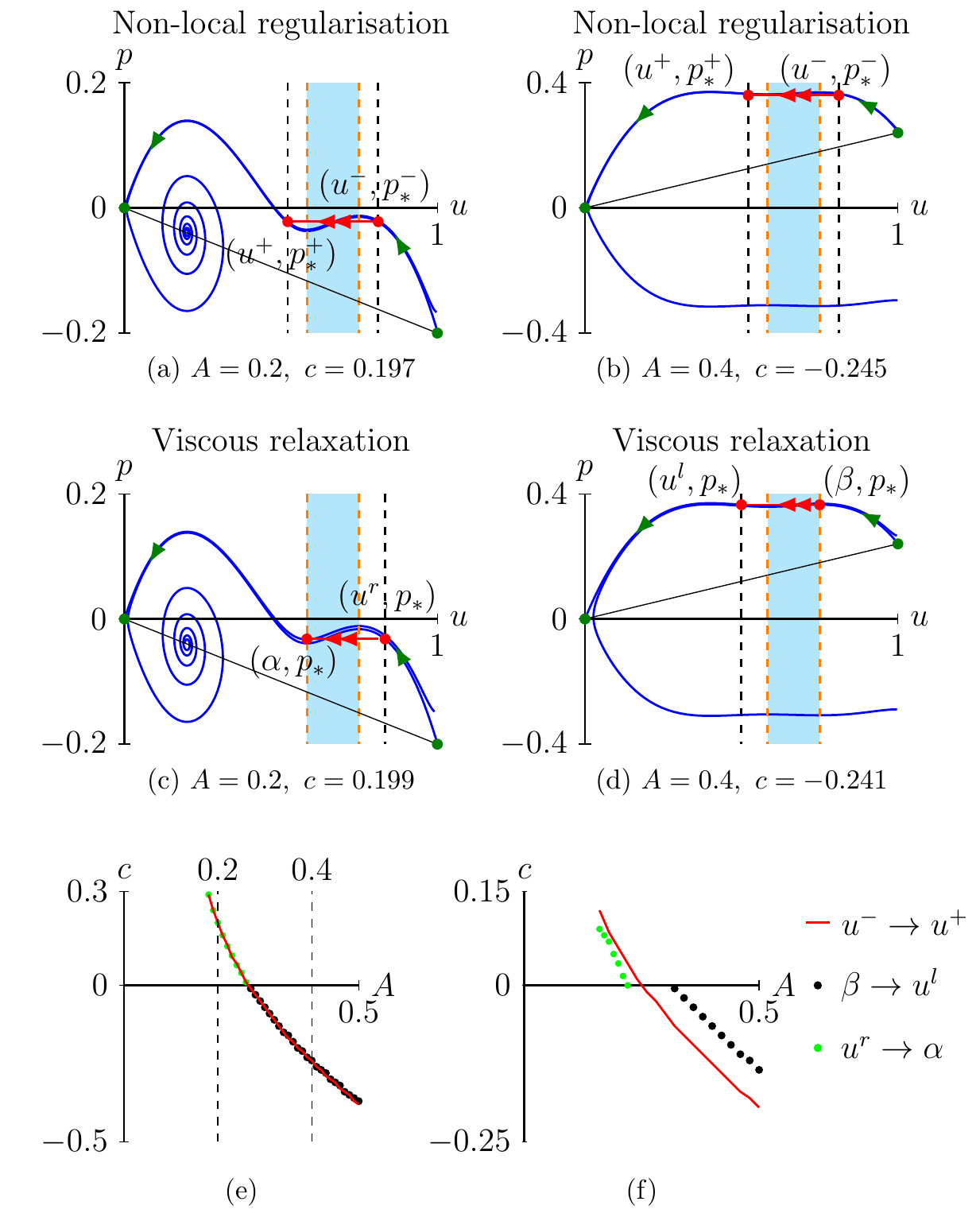}
\caption{(a) -- (d)
Phase planes of the desingularised system (\ref{GSPT_ALLEE_SYSTEMreduced problem3_2})/\eqref{GSPT_ALLEE_SYSTEMreduced problem3_3} with $D(u)=6(u-7/12)(u-3/4)$ and $R(u)=5u(1-u)(u-A)$ with $A=0.2$ (a) \& (c) and $A=0.4$ (b) \& (d). 
For $A=0.2$ in the non-local regularisation we observe a shock $u^-\to u^+$  for $c=0.197$, while for $A=0.2$ in the viscous relaxation we observe a shock $u^r\to\alpha$ for $c=0.199$. 
For $A=0.4$ in the non-local regularisation
we observe a shock $u^-\to u^+$ for $c=-0.245$, while for $A=0.4$ in the viscous relaxation we observe a shock $\beta\to u^l$ for $c=-0.241$. 
(e) The three different types of shock, $u^-\to u^+$, $u^r\to\alpha$ and $\beta\to u^l$, for changing $A$ from $0$ to $0.5$. The dashed line $A=0.2$ relates to (a) \& (b) and the dashed line $A=0.4$ relates to (c) \& (d). (f) For a $D(U)$ not centred around $2/3$ the differences among the speeds of the shocks is more prevalent. In particular, $D(u)=6(u-2/5)(u-3/4)$.}
\label{Fig.X}
\end{figure}

For the viscous relaxation the PDE \eqref{RDE_biology3} could be reduced to a singularly perturbed three-dimensional system of ODEs (\ref{Another4_1}). Whilst this ODE had the same reduced problem as with the non-local regularisation, it had a different layer problem \eqref{Another_layer_problem}. This difference can lead to shock-fronted travelling wave solutions with different characteristics for same nonlinear diffusivity function $D(U)$ (\ref{D(u)}) and reaction-term $R(U)$ (\ref{AlleeR(u)}), see Figure~\ref{Fig.X}. 
In addition, as the shock-connection in the layer problem is at a point where the critical manifold loses normal hyperbolicity, GSPT has to be extended to prove the persistence of the singular orbit for sufficiently small $\varepsilon$. Again, details of this computation were omitted.

\subsection{Regularisations and the lattice-based stochastic model}

While the two regularised PDEs have the same singular limit \eqref{RDE_biology1},  
the
different regularisations yielded shock-fronted travelling wave solutions with different characteristics.
Therefore, we mainly compared the singular limit results of the two models with the travelling wave ODE systems, and not with the numerical results of (\ref{RDE_biology1}).
The reason for this is that the numerical schemes used to simulate (\ref{RDE_biology1}) naturally introduce artificial regularisation (and error) terms and, as shown in this article, different regularisations yield shock-fronted travelling wave solutions with different characteristics. The connection between the numerical results of (\ref{RDE_biology1}) and the analytical results therefore needs to be further explored.
 
In addition, (\ref{RDE_biology1}) was derived from a lattice-based stochastic model and during this derivation of the continuous description small higher order terms were omitted. Including some of these small higher order terms
would potentially result in a (differently) regularised version of (\ref{RDE_biology1}), which in turn could lead to shock-fronted travelling wave solutions with different properties. 
Therefore, studying the connection between the lattice-based stochastic models and the regularisations is also an interesting topic.

For instance, a natural question to ask is what happens when we consider a linear combination of the
 non-local regularisation (considered in \S\ref{sec2}) and viscous regularisation (considered in \S\ref{sec3})
\begin{equation*}
     \begin{aligned}
    \frac{\partial U}{\partial t}=\frac{\partial}{\partial x}\left(D(U)\frac{\partial U}{\partial x}\right)+R(U)+(1-\mu)\varepsilon\frac{\partial^3U}{\partial x^2\partial t}-\mu\varepsilon^2\frac{\partial^4U}{\partial x^4},
    \end{aligned}
\end{equation*}
where $\mu\in[0,1]$ is a constant. Note that $\mu=0$ corresponds to the viscous regularisation (\ref{Another1}) and $\mu=1$ corresponds to the non-local regularisation~(\ref{RDE_biology3}). The associated four-dimensional slow-fast system\footnote{For $\mu=0$ this slow-fast system is actually three-dimensional and given by \eqref{Another4_1}.} is given by 
\begin{equation*}
    \left\{\begin{aligned}
    \varepsilon \mu &\frac{\textrm{d}u}{\textrm{d}z}&&=\left(
  \mu-1\right) cu+w,\\
  \varepsilon &\frac{\textrm{d}w}{\textrm{d}z}&&=v+F(u),\\
  &\frac{\textrm{d}p}{\textrm{d}z}&&=R(u),\\
  &\frac{\textrm{d}v}{\textrm{d}z}&&=p+cu.
\end{aligned}\right.
\end{equation*} 
The corresponding layer problem, for $\mu\neq 0$, is 
\begin{equation}
    \label{layer problem_two}
    \left\{\begin{aligned}
  &\frac{\textrm{d}u}{\textrm{d}\xi}&&=\left(1-\frac{1}{\mu}\right)cu+\frac{1}{\mu}w,\\
  &\frac{\textrm{d}w}{\textrm{d}\xi}&&=v+F(u).
\end{aligned}\right.
\end{equation}
If $v$ is such that \eqref{layer problem_two} has three fixed points $(u^-,w^-)$, $(u^0,w^0)$ and $(u^+,w^+)$, where $u^+<\alpha<u^0<\beta<u^-$. Then, for $\mu \neq 1$, \eqref{layer problem_two} does not have heteroclinic orbits connecting $(u^-,w^-)$ with $(u^+,w^+)$. Hence, we do not expect shock-fronted travelling wave solutions in this case.

\subsection{Generalisations}
\label{subsec2}
In this article, we concentrated on a specific quadratic nonlinear diffusivity function $D(U)$ (\ref{D(u)}) centred around $2/3$ and a specific Allee-type reaction-term $R(U)$ (\ref{AlleeR(u)}) as these were derived from an underlying lattice-based stochastic model \cite{johnston2017co}.
However, the techniques used in this article can in fact be easily extended to more general nonlinear diffusivity functions and reaction terms. For instance, if we change the reaction term from an Allee type (\ref{AlleeR(u)}) to a logistic type (\ref{logisticR(u)}) (as studied in \cite{2019arXiv190310090L}), we can still construct the higher-dimensional systems based on the two regularisations (\ref{RDE_biology3}) and (\ref{Another1}). Since the two layer problems \eqref{GSPT_ALLEE_SYSTEMlayer_problem} and \eqref{Another_layer_problem} only depend on $F(u)$, the anti-derivative of $D(u)$, and not on $R(u)$, we obtain the same conditions for the shocks as for the Allee type reaction term. That is, for the non-local regularisation the shocks will have, to leading order,  endpoints $u^-$ and $u^+$, while the shocks will have, to leading order,  endpoints $u^r$ and $\alpha$ or $u^l$ and $\beta$ for the viscous relaxation. In other words, the size of the shock depends on the relaxation and the nonlinear diffusivity function $D(U)$, but not the reaction term $R(U)$.
For both regularisations, the reduced desingularised problem has four fixed points which are determined by the roots of the product of the nonlinear diffusivity function $D(U)$ and the reaction term $R(U)$. In particular, the fixed points are $(0,0)$, $(1,-c)$, $(\alpha,-c\alpha)$ and $(\beta,-c\beta)$. In the desingularised system, the fixed point $(0,0)$ is a stable node or stable spiral for $c>0$ and an unstable node or unstable spiral for $c<0$. 
For shock-fronted travelling wave solutions with the asymptotic conditions $\lim_{z\to-\infty}U=1$ and $\lim_{z\to\infty}U=0$, we expect $(0,0)$ to be stable in the desingularised problem. Therefore, we expect those travelling wave solutions to have positive speeds. Hence, if the reaction term is logistic, we do not expect shock-fronted travelling wave solutions with negative speeds. However, using other boundary conditions may provide novel characteristics, see \cite{ELHACHEM2020132639,Fadai_2020} for examples of moving boundary problems with logistic type reaction terms.

\subsection{Stability}
\begin{figure}
\centering
\includegraphics[width=1\textwidth]{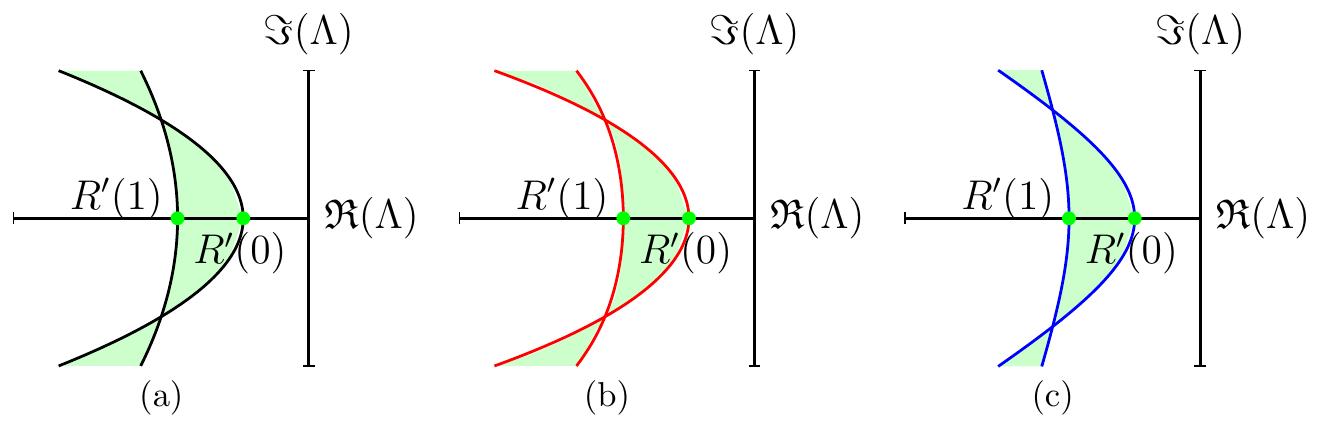}
\caption{The essential spectrums (shaded green regions plus boundaries) of the desingularised stability problems associated to \eqref{RDE_biology1} (a), \eqref{RDE_biology3} (b) and \eqref{Another1} (c) with $\varepsilon=0.1$. Note that the essential spectrums are to leading order the same and fully contained in the left-half plane.  
}
\label{esstre}
\end{figure}

Another natural extension of this work is to analyse the stability of the constructed shock-fronted travelling wave solutions. This was partly done for smooth travelling wave solution supported by \eqref{RDE_biology1} with $D(U)$ as in \eqref{D(u)} and logistic reaction term $R(U)$ \eqref{logisticR(u)} in \cite{2019arXiv190310090L}. In that article we studied the absolute spectrum of the associated desingularised stability problem and showed that for speeds above the minimal wave speed, the essential spectrum~\cite{kapitula2013spectral, sandstede2002stability} of the desingularised system can always be weighted into the left-half plane, while this is not possible for speeds below the minimal wave speed \cite{2019arXiv190310090L}. This analysis can be repeated for the shock-fronted travelling wave solutions constructed in this article since the essential spectrum is related to the behaviour of the wave at {\emph{infinity}} and thus only determined by the asymptotic end states of the shock-fronted travelling wave solution under consideration. For brevity we decided not to show this computation and instead refer to \cite{2019arXiv190310090L}. In short, the computation shows that the essential spectrum of the associated desingularised stability problems of (\ref{RDE_biology1}), (\ref{RDE_biology3}) and~(\ref{Another1}) are all fully contained in the left-half plane, see 
Figure \ref{esstre}, 
thus there are no absolute instabilities. However, what remains to be determined is the point spectrum, as well as the connection of the essential spectrum of the associated desingularised stability problem and the original stability problem, to complete the linear stability analysis. This is part of future work, see also the discussion in \cite{2019arXiv190310090L}.

\section*{Acknowledgements}
PvH, MW and MJS acknowledge support by the Australian Research Council (PvH: DP190102545 \& DP200102130, MJS: DP200100177, MW: DP200102130).

\vskip6pt

\enlargethispage{20pt}



\newpage
\appendix
\section{Necessary conditions for shock-fronted travelling wave solutions} \label{A:SH}
In this section, we follow \cite{Kuzmin2011819} and derive the two necessary conditions \eqref{necessarycondition1} and \eqref{necessarycondition2} for the existence of shock-fronted travelling wave solutions as mentioned in the Introduction. 
A shock-fronted travelling wave solutions of \eqref{RDE_biology1} solves the travelling wave ODE
\begin{equation*}
     \begin{aligned}
    c\frac{\textrm{d}u}{\textrm{d}z}+\frac{\textrm{d}}{\textrm{d} z}\left(D(u)\frac{\textrm{d}u}{\textrm{d}z}\right)+R(u)=0,
    \end{aligned}
\end{equation*}
where $z:=x-ct$ is the travelling wave coordinate. 
Define $g(u):=D(u)\textrm{d}u/\textrm{d}z$ in $(0,u_1)\cup(u_2,1)$, that is, $g(u)$ is defined in the region where the travelling wave solution $u$ is smooth. 
As we focus on monotonically decreasing travelling wave solutions we have that $g(u)<0$. The travelling wave ODE can now be written as 
\begin{equation}
    \label{AEQ1}
    g(u)\left(\frac{\textrm{d}(g(u))}{\textrm{d}u}+c\right)=-R(u)D(u).
\end{equation}
Integrating both sides \eqref{AEQ1} between $0$ and $u_a (< u_1)$ gives
\begin{equation}
\nonumber
    \int_{0}^{u_a} g(u)\textrm{d}g(u)+c\int_0^{u_a} g(u)\textrm{d}u=-\int_0^{u_a} R(u)D(u)\textrm{d}u,
\end{equation}
which leads to 
\begin{eqnarray*}
    c=-\dfrac{\displaystyle\int_0^{u_a} R(u)D(u)\textrm{d}u+\frac12(g(u_a))^2}{\displaystyle\int_0^{u_a} g(u)\textrm{d}u}.
\end{eqnarray*}
Thus, for $c<0$ a necessary condition for the existence of a shock-fronted travelling wave solution is
$$
\displaystyle\int_0^{u_a} R(u)D(u)\textrm{d}u<-\frac12(g(u_a))^2 < 0\,.
$$

Similarly, integrating \eqref{AEQ1} between $u_b (> u_2)$ and $1$ gives
\begin{equation*}
    c=-\dfrac{\displaystyle\int_{u_b}^1 R(u)D(u)\textrm{d}u-\frac12(g(u_b))^2}{\displaystyle\int_{u_b}^1 g(u)\textrm{d}u},
\end{equation*}
which, for $c>0$, leads to the necessary condition
$$
\displaystyle\int_{u_b}^1 R(u)D(u)\textrm{d}u>\frac12(g(u_b))^2>0
.
$$

\section{The heteroclinic orbits of the layer problem}
\label{A:HET}
We derive the analytic expressions for the heteroclinic orbits given in the layer problem supported by
\begin{equation*}
    \left\{\begin{aligned}
  &\frac{\textrm{d}u}{\textrm{d}\xi}&&=w,\\
  &\frac{\textrm{d}w}{\textrm{d}\xi}&&=v+F(u),
\end{aligned}\right.
\end{equation*}
where $v$ is a constant.
Based on its Hamiltonian, we require 
\begin{equation}
\nonumber
    H(u,w)=-\frac{1}{2}w^2+G(u)+vu=0,
\end{equation}
on the heteroclinic orbits $\left(u_h^{0,\pm},w_h^{0,\pm}\right)$. Subsequently, we obtain
\begin{equation}
\nonumber
    w=\pm\sqrt{2\left(G(u)+vu\right)}.
\end{equation}
Note that $G(u)$ has two integration constants. With specific integration constants, $w(u)$ can become a second-order polynomial with specific roots. That is, we can write $w$ as
\begin{equation}
    \nonumber
    w(u)=\pm\sqrt{2\left(G(u)+vu\right)}=\pm\sqrt{a^2(u-B_1)^2(u-B_2)^2}.
\end{equation}
Furthermore, as $w(u^{\pm})=0$, we can write $w$ as
\begin{equation}
    \nonumber
    w(u)=\pm a\left(u-u^+\right)\left(u-u^-\right),
\end{equation}
where $a=\sqrt{(D_i-D_g)/2}>0$. If we assume $w<0$ in $(u^+,u^-)$, then we have
\begin{equation}
\label{app_expression1}
    \frac{\textrm{d}u}{\textrm{d}\xi}=a(u-u^+)(u-u^-).
\end{equation}
Deriving the equation (\ref{app_expression1}) gives
\begin{equation}
\nonumber
    u_{h}^{0,+}(\xi)=\frac{u^++u^-}{2}+\frac{u^--u^+}{2}\tanh{\left(-\frac{a(u^--u^+)}{2}\xi\right)}.
\end{equation}
Subsequently, we obtain the expression of $w(\xi)$:
\begin{equation}
    \nonumber
    w_{h}^{0,+}(\xi)=-\frac{a(u^--u^+)^2}{4}\text{sech}^2{\left(-\frac{a(u^--u^+)}{2}\xi\right)},
\end{equation}
which satisfies $\lim_{\xi\to\pm\infty}w(\xi)=0$. Similarly, for the asymptotic conditions $\lim_{\xi\to-\infty}u(\xi)=u^+$ and $\lim_{\xi\to\infty}u(\xi)=u^-$, we have
\begin{equation}
     \label{app_expression2}
    \frac{\textrm{d}u}{\textrm{d}\xi}=-a(u-u^+)(u-u^-).
\end{equation}
Subsequently, solving (\ref{app_expression2}) gives the expressions of $u^{0,-}_h(\xi)$ and $w^{0,-}_h(\xi)$:
\begin{equation}
    \nonumber
    \begin{aligned}
    &u_{h}^{0,-}(\xi)=\frac{u^-+u^+}{2}-\frac{u^--u^+}{2}\tanh{\left(-\frac{a(u^--u^+)}{2}\xi\right)},\\
    &w_{h}^{0,-}(\xi)=\frac{a(u^--u^+)^2}{4}\text{sech}^2{\left(-\frac{a(u^--u^+)}{2}\xi\right)}.
    \end{aligned}
\end{equation}
\newpage

\bibliographystyle{elsarticle-harv} 
\bibliography{BIB_SHOCKS}

\end{document}